\documentclass[review,hidelinks,onefignum,onetabnum]{siamart220329}

\usepackage{lipsum}
\usepackage{amsfonts}
\usepackage{graphicx}
\usepackage{epstopdf}
\usepackage{algorithmic}
\ifpdf
  \DeclareGraphicsExtensions{.eps,.pdf,.png,.jpg}
\else
  \DeclareGraphicsExtensions{.eps}
\fi

\newsiamremark{remark}{Remark}
\newsiamremark{hypothesis}{Hypothesis}
\crefname{hypothesis}{Hypothesis}{Hypotheses}
\newsiamthm{claim}{Claim}

\headers{Least-squares neural network method}{Z. Cai, J. Choi, and M. Liu}

\title{Least-Squares Neural Network (LSNN) Method \\ for
Linear Advection-Reaction Equation: \\[1mm] Non-constant Jumps
\thanks{Submitted to the editors DATE.
\funding{This work was supported in part by the National Science Foundation under grant DMS-2110571.}}}

\author{Zhiqiang Cai\thanks{Department of Mathematics, Purdue University, 150 N. University Street, West Lafayette, IN 47907-2067 
  (\email{caiz@purdue.edu}, \email{choi508@purdue.edu}).}
\and Junpyo Choi\footnotemark[2]
\and Min Liu\thanks{School of Mechanical Engineering, Purdue University, 585 Purdue Mall,
West Lafayette, IN 47907-2088(\email{liu66@purdue.edu}). }}

\usepackage{amsopn}

\ifpdf
\hypersetup{
  pdftitle={Least-Squares Neural Network (LSNN) Method For
Linear Advection-Reaction Equation: Non-constant Jumps},
  pdfauthor={Z. Cai, J. Choi, and M. Liu}
}
\fi
\nolinenumbers
\usepackage[utf8]{inputenc}
\usepackage{bm}
\usepackage{mathtools}
\usepackage{indentfirst}
\usepackage{subfigure}
\usepackage{tcolorbox}
\usepackage{color}
\usepackage[export]{adjustbox}
\usepackage{scalerel}
\usepackage{tikz}
\usepackage{pgfplots}
\pgfplotsset{compat=1.18}
\usetikzlibrary{decorations.pathreplacing}

\usepackage{algorithmic}
\Crefname{ALC@unique}{Line}{Lines}

\newcommand{\R}{\mathbb{R}}

\usepackage{graphicx}
\usepackage{diagbox}
\usepackage{pifont}
\newcommand{\vertiii}[1]{{\left\vert\kern-0.25ex\left\vert\kern-0.25ex\left\vert #1 
    \right\vert\kern-0.25ex\right\vert\kern-0.25ex\right\vert}}

\usepackage{enumitem}
\setlist[itemize]{left=16pt} 

\def\bw{{\bf w}}
\def\bx{{\bf x}}

\def\cL{{\cal L}}
\def\cM{{\cal M}}

\def\cT{{\cal T}}

\begin{document}

\maketitle

\begin{abstract}
The least-squares ReLU neural network (LSNN) method was introduced and studied for solving linear advection-reaction equation with discontinuous solution in \cite{Cai2021linear,cai2023least}. The method is based on an equivalent least-squares formulation and \cite{cai2023least} employs ReLU neural network (NN) functions with $\lceil \log_2(d+1)\rceil+1$-layer representations for approximating solutions. In this paper, we show theoretically that the method is also capable of accurately approximating non-constant jumps along discontinuous interfaces that are not necessarily straight lines. Theoretical results are confirmed through multiple numerical examples with $d=2,3$ and various non-constant jumps and interface shapes, showing that the LSNN method with $\lceil \log_2(d+1)\rceil+1$ layers approximates solutions accurately with degrees of freedom less than that of mesh-based methods and without the common Gibbs phenomena along discontinuous interfaces having non-constant jumps.
\end{abstract}



\section{Introduction}
For decades, extensive research has been conducted on numerical methods for linear advection-reaction equations to develop precise and efficient numerical schemes. A major challenge in numerical simulation is that the solution of the equation is discontinuous along an interface because of a discontinuous inflow boundary condition, where the discontinuous interface can be the streamline from the inflow boundary. Traditional mesh-based numerical methods often exhibit oscillations near the discontinuity (called the Gibbs phenomena), which are not suitable for many applications and may not be extended to nonlinear hyperbolic conservation laws.

Recently, the application of neural networks (NNs) for solving partial differential equations have achieved significant accomplishments. For the linear advection-reaction problem, the least-squares ReLU neural network (LSNN) method was introduced and studied in \cite{Cai2021linear,cai2023least}. The method is based on an equivalent least-squares formulation studied in \cite{bochev2001improved,de2004least} and \cite{cai2023least} employs ReLU neural network (NN) functions with $\lceil \log_2(d+1)\rceil+1$-layer representations for approximating the solution. The LSNN method is capable of automatically approximating the discontinuous solution accurately since the free hyperplanes of ReLU NN functions adapt to the solution (see \cite{Cai2023nonlinear,Cai2021linear,cai2023least}). Moreover, for problems with unknown locations of discontinuity interfaces, it is quite easy to see that the LSNN method uses much fewer number of degrees of freedom than the mesh-based methods (see, e.g., \cite{Cai2023nonlinear,Cai2021linear}).

Approximation properties of ReLU NN functions to step functions were examined and employed in \cite{Cai2021linear, cai2023least}. In particular, it was shown theoretically that two- or $\lceil \log_2(d+1)\rceil+1$-layer ReLU NN functions are necessary and sufficient to approximate a step function with any given accuracy $\varepsilon>0$ when the discontinuous interface is a hyperplane or general hyper-surface, respectively. This approximation property was used to establish \textit{a priori} error estimates of the LSNN method for the linear advection-reaction problem.

The jump of the discontinuous solution of the problem, however, is generally  non-constant when the reaction coefficient is non-zero. The main purpose of this paper is to establish {\it a priori} error estimates (see \cref{main}) for the LSNN method without making the assumption that the jump is constant as in \cite{cai2023least}. To this end, we decompose the solution as the sum of discontinuous and continuous parts (see \cref{decop}), so that the discontinuous part of the solution can be described as a cylindrical surface on one subdomain and zero otherwise. Then we construct a continuous piecewise linear (CPWL) function with a sharp transition layer along the discontinuous interface to approximate the discontinuous piecewise cylindrical surface accurately. From \cite{tarela1999region,wang2005generalization,arora2016understanding,cai2023least}, we know that the CPWL function is a ReLU NN function $\mathbb{R}^d\to\mathbb{R}$ with a $\lceil \log_2(d+1)\rceil+1$-layer representation, from which it follows that the discontinuous part of the solution can be approximated by this class of functions for any prescribed accuracy. Then \cref{main} follows.

The rest of the paper is organized as follows. In \Cref{relu nn}, we introduce the linear advection-reaction problem, and briefly review and discuss properties of ReLU NN functions and the LSNN method in \cite{cai2023least}. Then theoretical convergence analysis is conducted in \Cref{error estimate}, showing that discretization error of the method for the problem mainly depends on the continuous part of the solution. Finally, to demonstrate the effectiveness of the method, we provide numerical results for test problems with various non-constant jumps in \Cref{numerics}. \Cref{conclusion} summarizes the work.

\section{Problem formulation and the LSNN method}\label{relu nn}
Let $\Omega$ be a bounded domain in ${\R}^d$ ($d\ge2$) 
with Lipschitz boundary $\partial \Omega$, and denote the advective velocity field by $\bm{\beta}(\bx) = (\beta_1, \cdots, \beta_d)^T \in C^0(\bar{\Omega})^d$. Define the inflow part of the boundary $\Gamma=\partial \Omega$ by
\begin{equation}
    \Gamma_- = \{\bx\in\Gamma :\, \bm{\beta}(\bx) \cdot \bm{n}(\bx) <0\}
\end{equation}
where $\bm{n}(\bx)$ is the unit outward normal vector to $\Gamma$ at $\bx\in \Gamma$. Consider the linear advection-reaction equation 
\begin{equation}\label{pde}
    \left\{\begin{array}{rccl}
    u_{\bm\beta} + \gamma\, u &=&f, &\text{ in }\, \Omega, \\[2mm]
    u&=&g, &\text{ on }\, \Gamma_{-},
    \end{array}\right.
\end{equation}
where $u_{\bm\beta}$ denotes the directional derivative of $u$ along $\bm{\beta}$. Assume that $\gamma \in C^0(\bar{\Omega})$, $f \in L^2(\Omega)$, and $g \in L^2(\Gamma_-)$ are given scalar-valued functions.

For the convenience of the reader, this section briefly reviews properties of ReLU neural network (NN) functions and the least-squares ReLU neural network (LSNN) method in \cite{cai2023least}. A function $\mathcal{N}:\mathbb{R}^d\to\mathbb{R}^c$ is called a ReLU neural network (NN) function if it can be expressed as a composition of functions
\begin{equation}\label{relu nn def}
 N^{(L)} \circ \cdots\circ N^{(2)}\circ N^{(1)}\text{ with }L>1,
\end{equation}
where $N^{(l)}:\mathbb{R}^{n_{l-1}}\to\mathbb{R}^{n_{l}}$ ($n_0=d$, $n_L=c$) is affine linear when $l=L$, and affine linear with the rectified linear unit (ReLU) activation function $\sigma$ applied to each component when $1\le l\le L-1$. Each affine linear function takes the form $\bm{\omega}^{(l)}\mathbf{x}-\mathbf{b}^{(l)}$ for $\mathbf{x}\in\mathbb{R}^{n_{l-1}}$ where $\bm{\omega}^{(l)}\in\mathbb{R}^{n_l\times n_{l-1}}$, $\mathbf{b}^{(l)}\in\mathbb{R}^{n_l}$ are weight and bias matrices, respectively. For $n\in\mathbb{N}$, denote the collection of all ReLU NN functions from $\mathbb{R}^d$ to $\mathbb{R}$ with depth $L$ and the number of hidden neurons $n(=n_1+\cdots+n_{L-1})$ by $\mathcal{M}(d,1,L,n)$ (1 being the output dimension), and the collection of all ReLU NN functions from $\mathbb{R}^d$ to $\mathbb{R}$ with depth $L$ by $\mathcal{M}(d,1,L)$. Then we have
\begin{equation}\label{relu dnn union}
  \mathcal{M}(d,1,L)=\bigcup_{n\in\mathbb{N}}\mathcal{M}(d,1,L,n).  
\end{equation}

The following proposition justifies our use of $\lceil \log_2(d+1)\rceil+1$-layer ReLU NN functions.

\begin{proposition}[see \cite{arora2016understanding,cai2023least}]\label{cpwl=relu}
The collection of all continuous piecewise linear {\em(CPWL)} functions on $\mathbb{R}^d$ is equal to $\cM(d,1,\lceil \log_2(d+1)\rceil+1)$, i.e., the collection of all {\em ReLU NN} functions from $\mathbb{R}^d$ to $\mathbb{R}$ that have representations with depth $\lceil \log_2(d+1)\rceil+1$.
\end{proposition}

\begin{proposition}\label{subset}
$\cM(d,1,L,n)\subseteq\cM(d,1,L,n+1)$.
\end{proposition}

Proposition \ref{subset} implies that as we increase $n$, $\cM(d,1,L,n)$ approaches $\mathcal{M}(d,1,L)$ and the approximation class gets larger.

As in \cite{cai2023least}, breaking hyperplanes are depicted in \cref{test1 figure,test2 figure,test3 figure,test4 figure,test5 figure,test6 figure} to better understand the graphs of ReLU NN function approximations using domain partitions (on each element in a given partition, the ReLU NN function is affine linear; see \cite{cai2023least}). More specifically, the $l^{\text{th}}$- (hidden) layer breaking hyperplanes of a given ReLU NN function (with output dimension 1) representation as in \eqref{relu nn def} are defined as the zero sets of the layer without the activation function: when $l=1$, $\mathbf{w}_i^{(1)}\mathbf{x}-b_i^{(1)}=0$, and when $2\le l<L$,
\[\mathbf{w}_i^{(l)}\left(N^{(l-1)}\circ\cdots\circ N^{(2)}\circ N^{(1)}(\mathbf{x})\right)-b_i^{(l)}=0,\]

where
\[\bm{\omega}^{(l)}=(\mathbf{w}_1^{(l)},\ldots,\mathbf{w}_{n_l}^{(l)})^T\in\mathbb{R}^{n_l\times n_{l-1}},\quad\text{and}\quad \mathbf{b}^{(l)}=(b_1^{(l)},\ldots,b_{n_l}^{(l)})^T.\]

Define the least-squares (LS) functional 
 \begin{equation}\label{ls}
    \mathcal{L}(v;{\bf f}) = \|v_{\bm\beta} +\gamma\, v-f\|_{0,\Omega}^2 +  \|v-g\|_{-\bm\beta}^2, 
\end{equation}
where ${\bf f} = (f,g)$ and $\|\cdot\|_{-\bm\beta}$ is given by
 \[
 \|v\|_{-\bm{\beta}} 
 =\left<v,v\right>^{1/2}_{-\bm{\beta}} 
 =\left( \int_{\Gamma_-} |\bm{\beta}\! \cdot \!\bm{n}|\, v^2\,ds\right)^{1/2}.
 \]
The LS formulation of problem \cref{pde} is to seek $u\in V_{\bm\beta}$ such that
\begin{equation}\label{minimization1}
    \mathcal{L}(u;{\bf f}) = \min_{\small v\in V_{\bm\beta}} \mathcal{L}(v;{\bf f}),
\end{equation}
where $V_{\bm\beta} = \{v\in L^2(\Omega): v_{\bm{\beta}}\in L^2(\Omega)\}$ is a Hilbert space that is equipped with the norm
 \[
 \vertiii{v}_{\bm\beta}= \left(\|v\|_{0,\Omega}^2 + \|v_{\bm\beta}
 \|_{0,\Omega}^2 \right)^{1/2}.
 \]
Then the corresponding LS and discrete LS approximations are, respectively, to find $u_{_N} \in \cM(d,1,L,n)$ such that 
\begin{equation}\label{L-NN}
     \mathcal{L}\big(u_{_N};{\bf f}\big)
     = \min\limits_{v\in \cM(d,1,L,n)} \mathcal{L}\big(v;{\bf f}\big),
\end{equation}
and to find ${u}^{_N}_{_{\small {\cal T}}}\in \cM(d,1,L,n)$ such that
 \begin{equation}\label{discrete_minimization_functional}
  \mathcal{L}_{_{\small {\cal T}}} \big({u}^{_N}_{_{\small {\cal T}}};{\bf f}\big) 
  = \min\limits_{v\in \cM(d,1,L,n)} \mathcal{L}_{_{\small {\cal T}}}\big(v;{\bf f}\big),
\end{equation}
where $\mathcal{L}_{_{\small {\cal T}}}\big(v;{\bf f}\big)$ is the discrete LS functional (see \cite{Cai2021linear,cai2023least}).

\section{Error estimates}\label{error estimate}
In this section, we establish error estimates of the LSNN method for the linear advection-reaction equation with a non-constant jump along a discontinuous interface. For simplicity, we will restrict our attention to two dimensions.

\begin{figure}[htbp]
\centering
\subfigure[Subdomains $\Upsilon_1$, $\Upsilon_2$ with respect to $\bm\beta$\label{subdomains}]{
\begin{minipage}[t]{0.4\linewidth}
\centering
\begin{tikzpicture}[scale=0.8, transform shape]
    \draw [] (0,5)-- (5,5);
    \draw [] (0,5)-- (0,0);
    \draw [] (0,0)-- (5,0);
    \draw [] (5,0)-- (5,5);
    \draw [] (0,3.5)-- (2.5,2.5);
    \draw [] (3.5,0)--(2.5,2.5);
    \draw [dotted, very thick] (0,0)--(5,5);
    
    \node[] at (-0.4,0) {$\Omega_{1}$};
    \node[] at (5.4,5) {$\Omega_{2}$};
    \node[] at (2.5,1.9) {$I$};
    \node[] at (2,4) {$\Upsilon_1$};
    \node[] at (4,2) {$\Upsilon_2$};
    \node[] at (3.5,-0.4) {$\mathbf{x}_0$};
    \filldraw[black] (3.5,0) circle (2pt);
    \end{tikzpicture}
\end{minipage}%
}%
\hspace{0.2in}
\subfigure[An $\varepsilon_2$ neighborhood around the interface $I$ in the direction of $\bm\beta$\label{nbd}]{
\begin{minipage}[t]{0.4\linewidth}
\centering
\begin{tikzpicture}[scale=0.8, transform shape]
    \draw [] (0,5)-- (5,5);
    \draw [] (0,5)-- (0,0);
    \draw [] (0,0)-- (5,0);
    \draw [] (5,0)-- (5,5);
    \draw [] (0,3.5)-- (2.5,2.5);
    \draw [] (3.5,0)--(2.5,2.5);
    \draw [dotted, very thick] (0,0)--(5,5);
    \draw [dotted, very thick] (15/7,15/7)--(3,0);
    \draw [dotted, very thick] (15/7,15/7)--(0,3);
    \draw [dotted, very thick] (20/7,20/7)--(4,0);
    \draw [dotted, very thick] (20/7,20/7)--(0,4);
    \draw [] (18/7,1.07)--(3,1.25);
    \draw[-latex] (1.9,0.5).. controls (2.6,0.4) ..(2.84,1.1);

    \node[] at (1.7,0.6) {$\varepsilon_2$};
    \node[] at (-0.4,0) {$\Omega_{1}$};
    \node[] at (5.4,5) {$\Omega_{2}$};
    \node[] at (2.5,1.9) {$I$};
    \node[] at (2,4) {$\Upsilon_1$};
    \node[] at (4,2) {$\Upsilon_2$};
    \end{tikzpicture}
\end{minipage}%
}%
\\
\subfigure[A subdomain $\Upsilon_i$\label{one sub}]{
\begin{minipage}[t]{0.4\linewidth}
\centering
\begin{tikzpicture}[scale=0.8, transform shape]
    \draw [] (0,5)-- (5,5);
    \draw [] (0,5)-- (0,0);
    \draw [] (0,0)-- (5,0);
    \draw [] (5,0)-- (5,5);
    \draw [] (2.5,0)--(2.5,5);
    \draw [] (2.5,1.5)--(3,1.5);
    \draw [dotted, very thick] (2,0)--(2,5);
    \draw [dotted, very thick] (3,0)--(3,5);
    \draw[-latex] (3.5,1).. controls (3,0.5) ..(2.75,1.4);

    \node[] at (-0.4,0) {$\Upsilon_i$};
    \node[] at (2.25,2.5) {$I$};
    \node[] at (1,2.5) {$\Omega_{1i}$};
    \node[] at (4,2.5) {$\Omega_{2i}$};
    \node[] at (3.7,1.2) {$\varepsilon_2$};
    \node[] at (2.5,-0.4) {$(x_0,0)$};
    \filldraw[black] (2.5,0) circle (2pt);
    \end{tikzpicture}
\end{minipage}%
}%
\caption{A domain decomposition for the case that $\bm\beta$ is piecewise constant}
\end{figure}
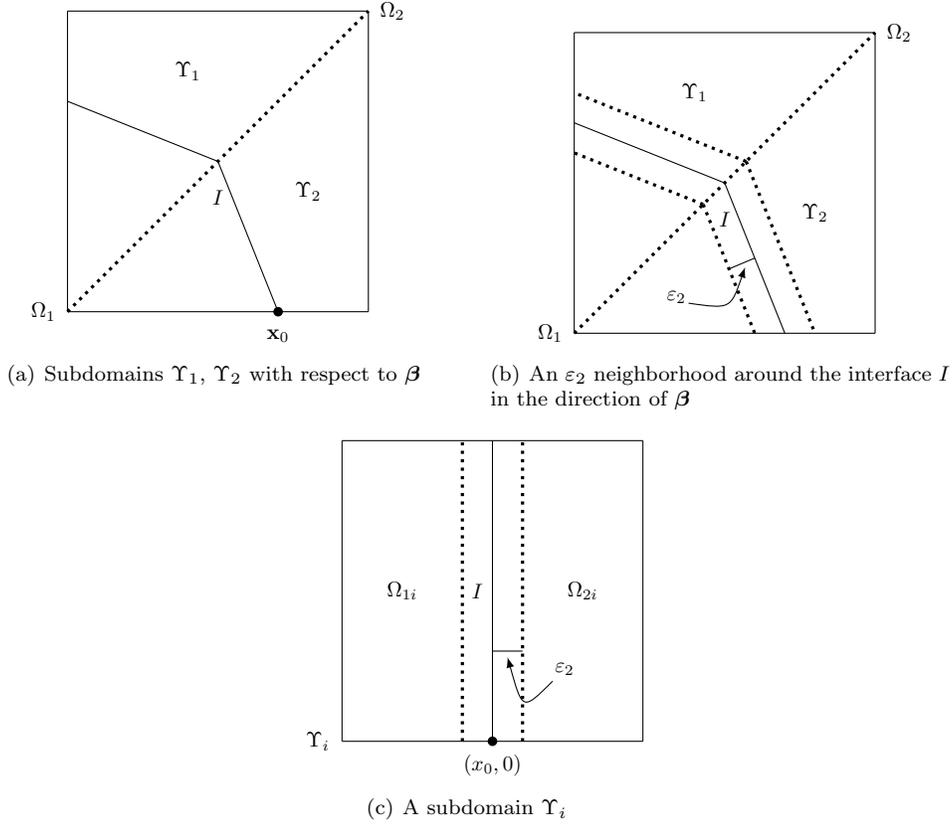

To this end, assume the advection velocity field $\bm\beta$ is piecewise constant. That is, there exists a partition of the domain $\Omega$ such that $\bm\beta$ has the same direction but possibly a different magnitude at each interior point of each subdomain. Without loss of generality, assume that there are only two sub-domains: $\Omega=\Upsilon_1\cup\Upsilon_2$ and that the inflow boundary data $g(\bx)$ is discontinuous at only one point $\bx_0\in \Gamma_-$ with $g(\bx^-_0)=\alpha_1$ and $g(\bx^+_0)=\alpha_2$ from different sides. (\cref{subdomains} depicts $\Upsilon_1$ and $\Upsilon_2$ as the left-upper and the right-lower triangles, respectively.) 

Let $I$ be the streamline emanating from $\bx_0$; then the discontinuous interface $I$ divides the domain $\Omega$ into two sub-domains: $\Omega=\Omega_1\cup \Omega_2$, where $\Omega_1$ and $\Omega_2$ are the left-lower and the right-upper subdomains separated by the discontinuous interface $I$, respectively (see \cref{subdomains}). The corresponding solution $u$ of \eqref{pde} is discontinuous across the interface $I$ and is piecewise smooth with respect to the partition $\{\Omega_1,\Omega_2\}$.
For a given $\varepsilon_2>0$, take an $\varepsilon_2$ neighborhood around the interface $I$ in the direction of $\bm\beta$ as in \cref{nbd}.

Next, we estimate the error in the sub-domain $\Upsilon_i$, say, $\Upsilon_2$.
To further simplify the error estimate, we assume that 
\[
\Upsilon_i=(0,1)\times(0,1), \quad 
\bm\beta(\mathbf{x}) = (0,v_2(\mathbf{x}))^T, \quad\mbox{and}\quad \mathbf{x}_0=(x_0,0) \mbox{ for } x_0\in (0,1).
\]
These assumptions imply that the restriction of the interface $I$ to $\Upsilon_i$ is a vertical line segment
\[
I=\{(x_0,y)\in\Upsilon_i:y\in(0,1)\},
\]
and that $\Upsilon_i$ is partitioned into two sub-domains
\[
\Omega_{1i}=\{\bx=(x,y)\in\Upsilon_i :\, x <x_0\}
\quad\mbox{and}\quad
\Omega_{2i}=\{\bx=(x,y)\in\Upsilon_i :\, x >x_0\}
\]
(see \cref{one sub}).

In $\Upsilon_i=\Omega_{1i}\cup \Omega_{2i}\cup I$, let $u_1$ and $u_2$ be the solutions of \cref{pde} defined only on $\Omega_{1i}$ with the constant inflow boundary conditions ${g}=\alpha_1$ and $\alpha_2$ on $\{(x,0) : \, x\in [0,x_0]\}$, respectively. (When $\Upsilon_i$ is different from $\Upsilon_2$, the discontinuous point is not $\mathbf{x}_0$ but an interior point of the domain $\Omega$, and the values of the solution $u$ at that discontinuous point from different sides are taken as the constant inflow boundary conditions.) 
We set $a(\bx)=u_1(\bx)-u_2(\bx)$ and let $\chi(\bx)$ be the piecewise discontinuous function defined by 
\begin{equation}\label{chi}
 \chi(\bx)=\left\{\begin{array}{rl}
 a(\bx), & \bx \in \Omega_{1i},\\[2mm]
 0, & \bx \in \Omega_{2i};
 \end{array}
 \right.   
\end{equation}
then the solution $u$ of \cref{pde} has the following decomposition (see \cite{cai2023least})
\begin{equation}\label{decop}
u(\bx)=\hat{u} (\bx) + \, \chi (\bx)\,\text{ in }\Upsilon_i.   
\end{equation}
Here, $\hat{u}(\bx)=u(\bx)- \, \chi (\bx)$ is clearly piecewise smooth; moreover, it is also continuous in $\Upsilon_i$ since $\hat{u}\big|_{I}= u_2\big|_{I}$ from both sides. Then we have the following error estimate and postpone its proof to Appendix.
\begin{theorem}\label{thm:chi}
For any $\varepsilon_2>0$ and $\varepsilon_3>0$, on $\Upsilon_i$, there exists a CPWL function $p_i(\bx)$ such that
\begin{equation}
    \vertiii{\chi-p_i}_{\bm\beta}\le D_1\sqrt{\varepsilon_2}+D_2\sqrt{\varepsilon_3}.
\end{equation}
\end{theorem}

\begin{remark}
We now construct the CPWL function $p(\mathbf{x})$ on $\Omega$ defined by
\[p(\mathbf{x})=p_i(\mathbf{x}),\ \mathbf{x}\in\Upsilon_i,\]
such that $p_i(\mathbf{x})=p_{i+1}(\mathbf{x})$ on the intersection of $\Upsilon_i$ and $\Upsilon_{i+1}$.
Using the triangle inequality, \cref{thm:chi} can be extended to the case that $\bm\beta$ is piecewise constant to establish the error estimate on the whole domain $\Omega$.
\end{remark}

\begin{theorem}\label{main}
Let $u$ and $u_{_{N}}$ be the solutions of problems \cref{minimization1} and \cref{L-NN}, respectively. If the depth of ReLU NN functions in \cref{L-NN} is at least $\lceil \log_2(d+1)\rceil+1$, then for a sufficiently large integer $n$, there exists an integer $\hat{n}\leq n$ such that
\begin{equation}\label{tau_1-error-2}
 \vertiii{u-u_{_{N}}}_{\bm\beta}
 \leq C\,\left(\sqrt{\varepsilon_2}+\sqrt{\varepsilon_3} + \inf_{v\in \cM(d,n-\hat{n})} \vertiii{\hat{u}-v}_{\bm\beta}
 \right),
 \end{equation}
 where $\cM(d,n-\hat{n})=\cM(d,1,\lceil \log_2(d+1)\rceil+1,n-\hat{n})$.
\end{theorem}

\begin{proof}
The proof is similar to that of Theorem 4.4 in \cite{cai2023least}.
\end{proof}

\begin{lemma}\label{uutn}
Let $u$, $u_{_N}$, and $u_{_\cT}^{_N}$ be the solutions of problems \cref{minimization1}, \cref{L-NN}, and \cref{discrete_minimization_functional}, respectively. Then there exist positive constants $C_{1}$ and $C_{2}$ such that
\begin{equation}\label{Cea-L-d}
\begin{split}
     \vertiii{u-u^{_N}_{_\cT}}_{\bm\beta}
    \le\  &C_{1}\,\left(\big|(\cL-\cL_{_\cT})(u_{_N}-u_{_\cT}^{_N}, {\bf 0})\big|
    + \big|(\cL-\cL_{_\cT})(u-u_{_N}, {\bf 0})\big|
    \right)^{1/2}\\
    &+C_{2}\,\left(\sqrt{\varepsilon_2}+\sqrt{\varepsilon_3}+ \inf_{v\in \cM(d,n-\hat{n})} \vertiii{\hat{u}-v}_{\bm\beta}
 \right).
\end{split}
\end{equation}
\end{lemma}

\begin{proof}
The proof is similar to that of Lemma 4.7 in \cite{cai2023least}.
\end{proof}

\section{Numerical experiments}\label{numerics}
In this section, we demonstrate the performance of the LSNN
method across different settings, incorporating various non-constant jumps. The discrete LS functional was minimized by the Adam optimization algorithm \cite{kingma2015} on a uniform mesh with mesh size $h=10^{-2}$. As in \cite{cai2023least}, the directional derivative $v_{\bm\beta}$ was approximated by the backward finite difference quotient multiplied by $|{\bm\beta}|$
\begin{equation}\label{finite_diff}
    v_{\bm\beta}(\bx_{_K}) \approx|{\bm\beta}| \frac{v(\bx_{_K})-v\big(\bx_{_K} - \rho\bar{\bm{\beta}}(\bx_{_K})\big)}{\rho},
\end{equation}
where $\bar{\bm\beta}=\frac{\bm\beta}{|\bm\beta|}$ and $\rho=h/4$ (except for the fifth test problem, which used $\rho=h/15$). The LSNN method was implemented with an adaptive learning rate that started with $0.004$ and was reduced by half for every 50000 iterations (except for the fourth test problem, which reduced for every 100000 iterations). For each experiment, to avoid local minima, 10 ReLU NN functions were trained for 5000 iterations each, and then the experiment began with one of the pretrained network functions that gave the minimum loss. 

\cref{test1 table,test2 table,test3 table,test4 table,test5 table,test6 table} report numerical errors in the relative $L^2$, $V_{\bm\beta}$, and the LS functional with parameters being the total number of weights and biases. Since the input dimensions $d=2,3$ and the depth $\lceil \log_2(d+1)\rceil+1=3$ for $d=2,3$, we employed ReLU NN functions with 2--$n_1$--$n_2$--1 or 3--$n_1$--$n_2$--1 representations or structures, which means that the representations have two-hidden-layers with $n_1$, $n_2$ neurons, respectively (here 2,3 mean the input dimensions and 1 is the output dimension.) In the fourth test problem for which the discontinuous interface is not a straight line, we also have the approximation of the 2-layer NN known as a universal approximator (see, e.g., \cite{leshno1993multilayer,pinkus1999approximation}) to show how the depth of a neural network impacts the approximation (see \cite{cai2023least} for more examples).

All the test problems are defined on the domain $\Omega=(0,1)^2$ or $(0,1)^3$ with $\gamma=1$ ($f=1$ for the first three and the last test problems, and $f=0$ for the remaining test problems).

\subsection{A problem with a constant advection velocity field}\label{test1}
The first test problem has the constant advective velocity field $\bm{\beta}(x,y) = (0,1),\,\, (x,y)\in\Omega$, and the inflow boundary of the problem is $\Gamma_{-}=\{(x,0):x\in(0,1)\}$.
The inflow boundary condition is given by
\begin{equation*}
g(x,y)=\left\{ \begin{array}{rl}
 1,& (x,y)\in \Gamma^1_-\equiv \{(x,0): x\in(0,1/2)\}, \\[2mm]
 2, &(x,y)\in \Gamma^2_-=\Gamma_-\setminus \Gamma_-^1.
\end{array}\right.
\end{equation*} 
The exact solution of this test problem is
\begin{equation*}
u(x,y)=\left\{ \begin{array}{rl}
 1,& (x,y)\in \Omega_1=\{(x,y)\in\Omega:x< 1/2\}, \\[2mm]
 1+e^{-y}, & (x,y)\in \Omega\setminus\Omega_1.
\end{array}\right.
\end{equation*}

\begin{figure}[htbp]\label{test1 figure}
\centering
\subfigure[The interface\label{interface1}]{
\begin{minipage}[t]{0.4\linewidth}
\centering
\includegraphics[width=1.8in]{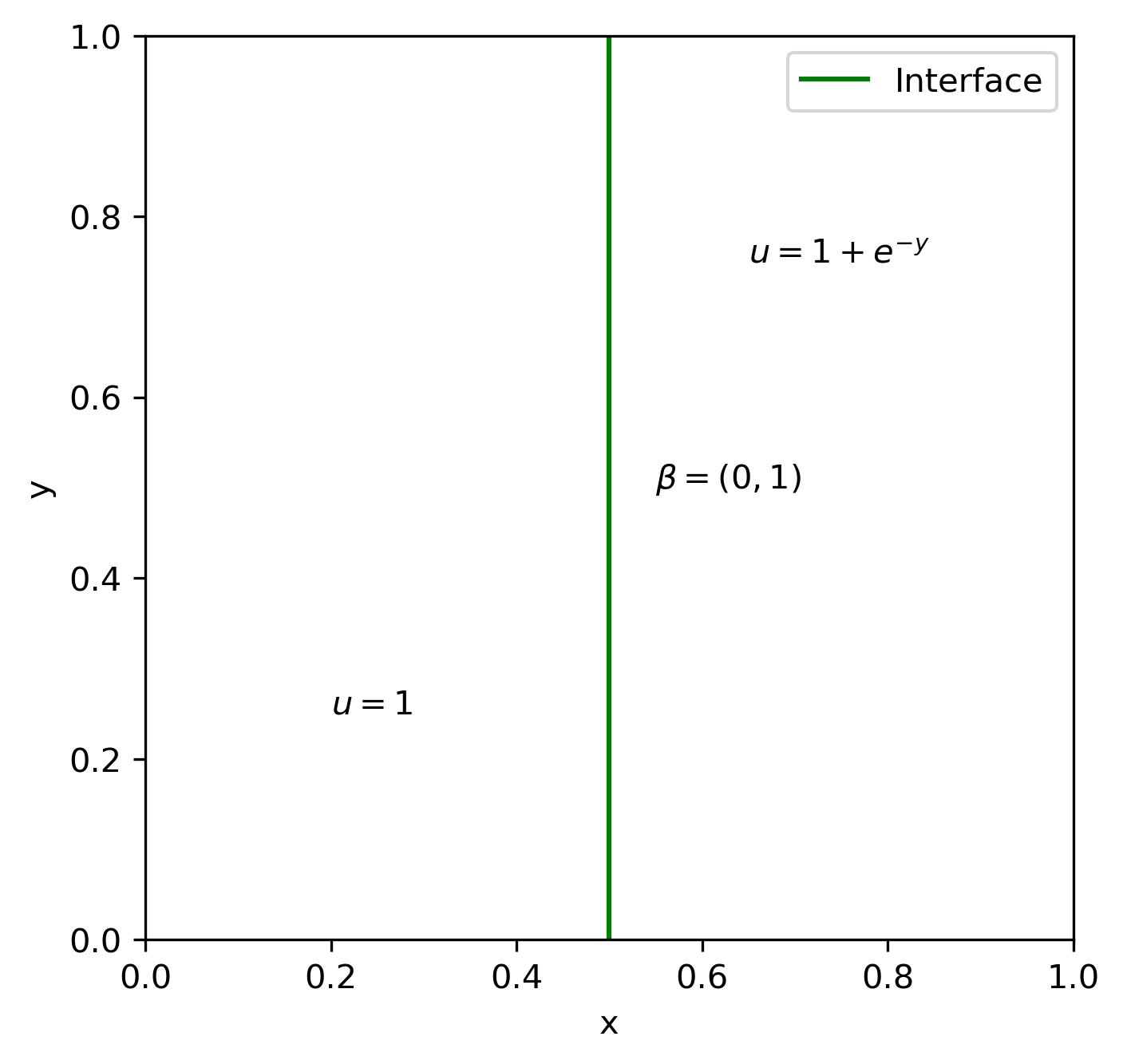}
\end{minipage}%
}%
\hspace{0.2in}
\subfigure[The trace of Figure \ref{comparison1} on $y=0.5$\label{vertical1}]{
\begin{minipage}[t]{0.4\linewidth}
\centering
\includegraphics[width=1.8in]{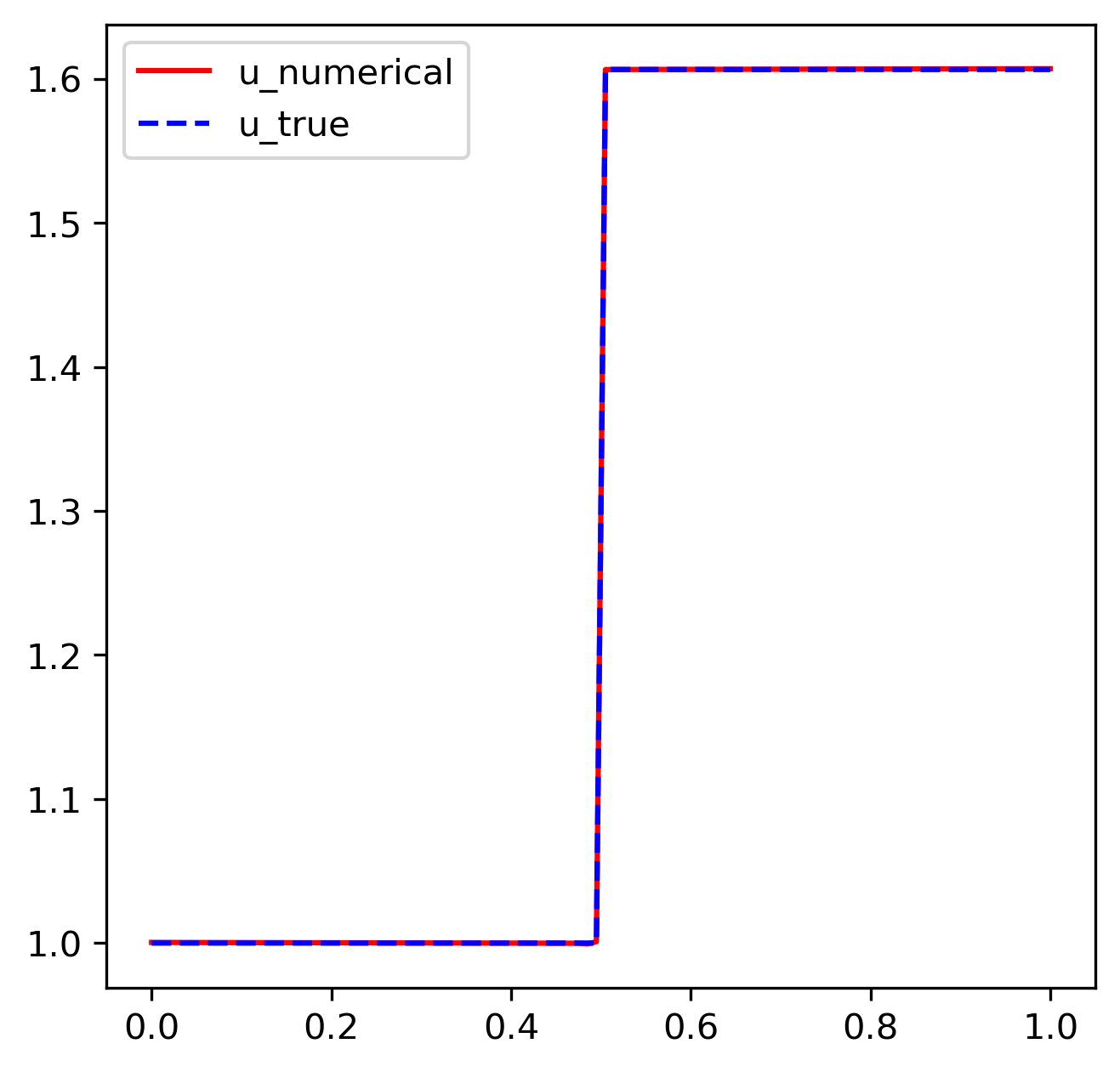}
\end{minipage}%
}%
\\
\subfigure[The exact solution\label{comparison_exact1}]{
\begin{minipage}[t]{0.4\linewidth}
\centering
\includegraphics[width=1.8in]{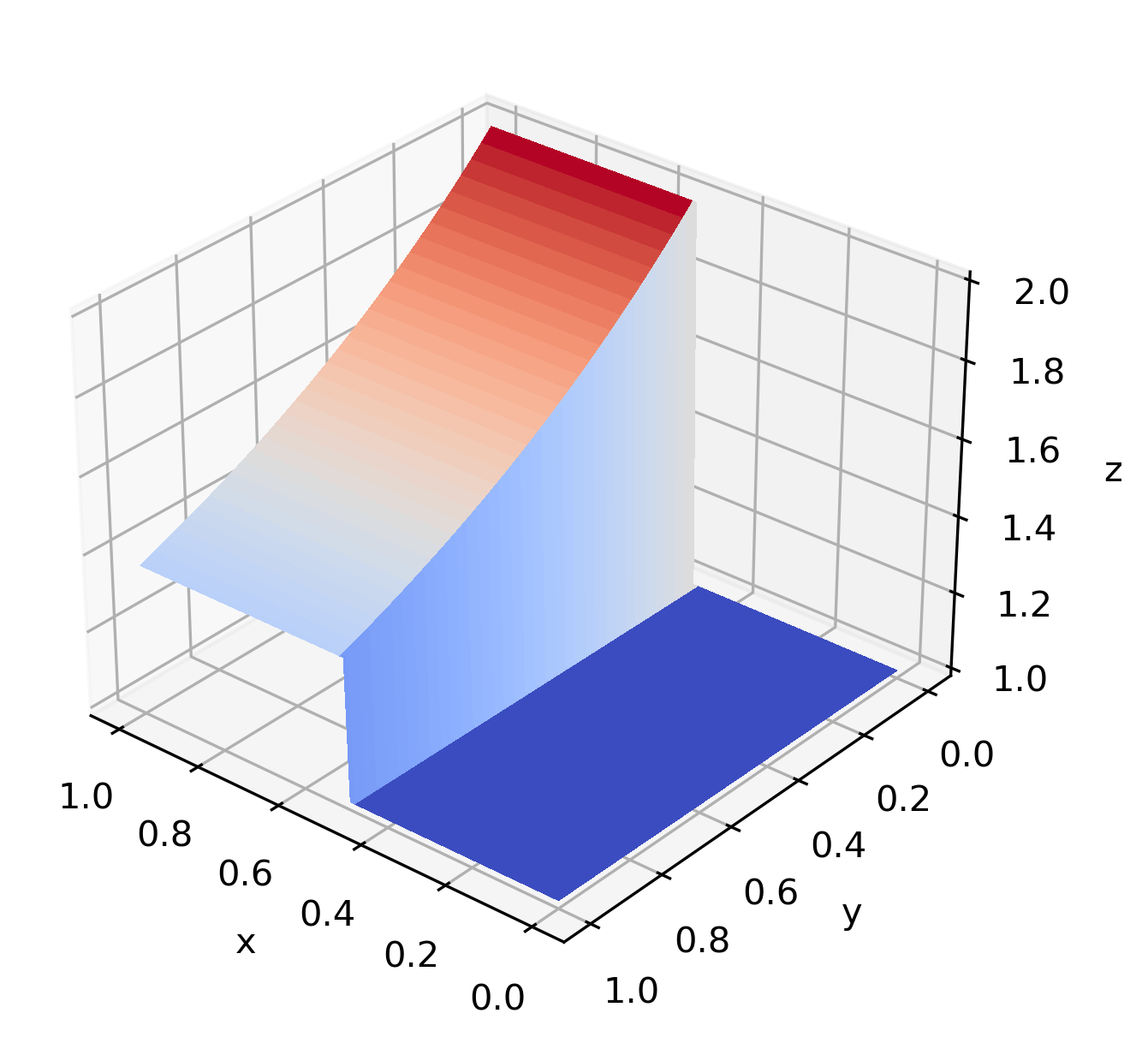}
\end{minipage}%
}%
\hspace{0.2in}
\subfigure[A 2--20--20--1 ReLU NN function approximation\label{comparison1}]{
\begin{minipage}[t]{0.4\linewidth}
\centering
\includegraphics[width=1.8in]{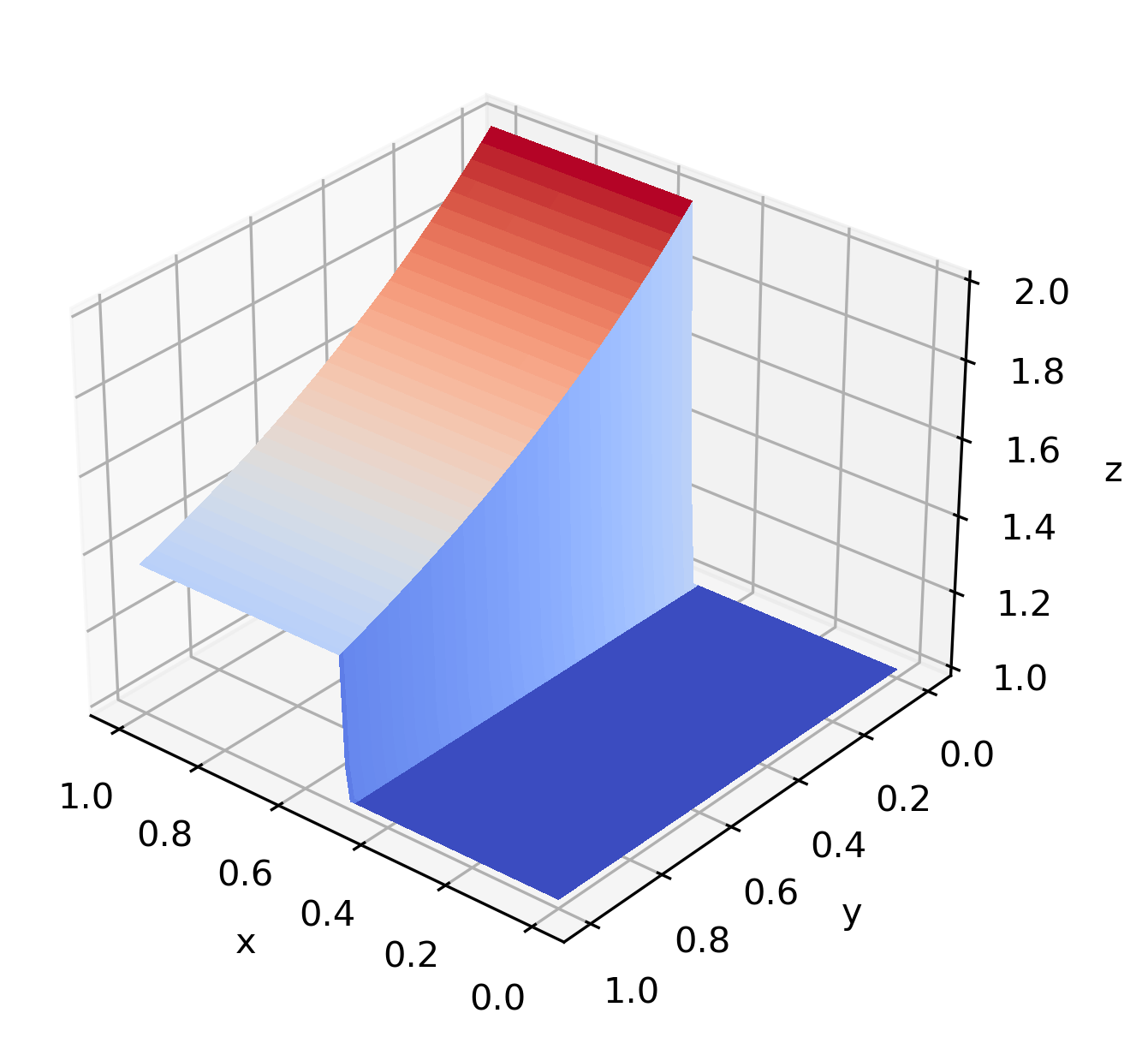}
\end{minipage}%
}%
\\
\subfigure[The breaking hyperplanes of the approximation in Figure \ref{comparison1}\label{breaking1}]{
\begin{minipage}[t]{0.4\linewidth}
\centering
\includegraphics[width=1.8in]{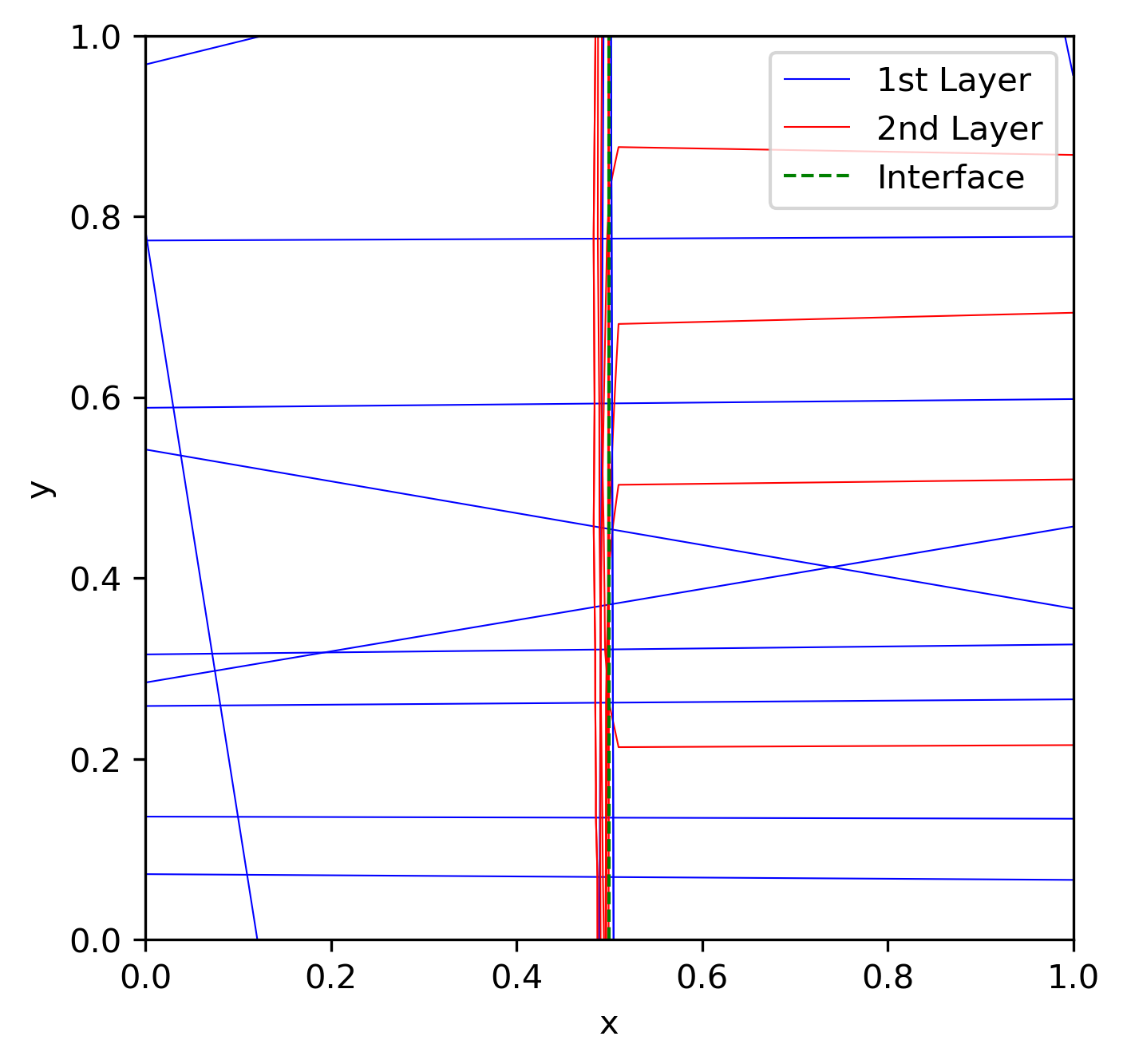}
\end{minipage}%
}%
\caption{Approximation results of the problem in \Cref{test1}}
\end{figure}

The LSNN method was implemented with 50000 iterations for 2--20--20--1 ReLU NN functions. We report the numerical results in \cref{test1 figure,test1 table}. The traces (\cref{vertical1}) of the exact and numerical solutions on the plane $y=0.5$ show no difference or oscillation. The exact solution (\cref{comparison_exact1}), which has a non-constant jump along the vertical interface (\cref{interface1}) is accurately approximated by a 3-layer ReLU NN function (\cref{comparison1,test1 table}). We note that the solution of this test problem takes the same form as $\chi(\mathbf{x})$ in \eqref{chi}, which was approximated by a CPWL function constructed by partitioning the domain into rectangles stacking on top of each other. It appears from \cref{breaking1} that the 3-layer ReLU NN function approximation has a similar partition, and the second-layer breaking hyperplanes were generated for approximating the jump along the discontinuous interface and the non-constant part of the solution, which is consistent with our theoretical analysis on the convergence of the method.

\begin{table}[htbp]\label{test1 table}
\caption{Relative errors of the problem in \Cref{test1}}
\centering
\begin{tabular}{|l|l|l|l|l|}
\hline
Network structure  &$\frac{\|u-{u}^{_N}_{_\cT}\|_0}{\|u\|_0}$ &$\frac{\vertiii{u-{u}^{_N}_{_\cT}}_{\bm\beta}}{\vertiii{u}_{\bm\beta}}$ & $\frac{\mathcal{L}^{1/2}({u}^{_N}_{_\cT},\bf f)}{\mathcal{L}^{1/2}({u}^{_N}_{_\cT},\bf 0)}$ & Parameters \\ \hline
2--20--20--1  & 0.037881 & 0.007044 & 0.005391   & 501\\ \hline
\end{tabular}
\end{table}

\subsection{A problem with a piecewise smooth solution}\label{test2}
This example is a modification of \Cref{test1} by changing the inflow boundary condition to
\begin{equation*}
g(x,y)=\left\{ \begin{array}{rl}
 0,& (x,y)\in \Gamma^1_-\equiv \{(x,0): x\in(0,1/2)\}, \\[2mm]
 2, &(x,y)\in \Gamma^2_-=\Gamma_-\setminus \Gamma_-^1.
\end{array}\right.
\end{equation*} 
The exact solution of this test problem is
\begin{equation}
u(x,y)=\left\{ \begin{array}{rl}
 1-e^{-y},& (x,y)\in \Omega_1, \\[2mm]
 1+e^{-y}, & (x,y)\in \Omega_2.
\end{array}\right.
\end{equation}

\begin{figure}[htbp]\label{test2 figure}
\centering
\subfigure[The interface\label{interface2}]{
\begin{minipage}[t]{0.4\linewidth}
\centering
\includegraphics[width=1.8in]{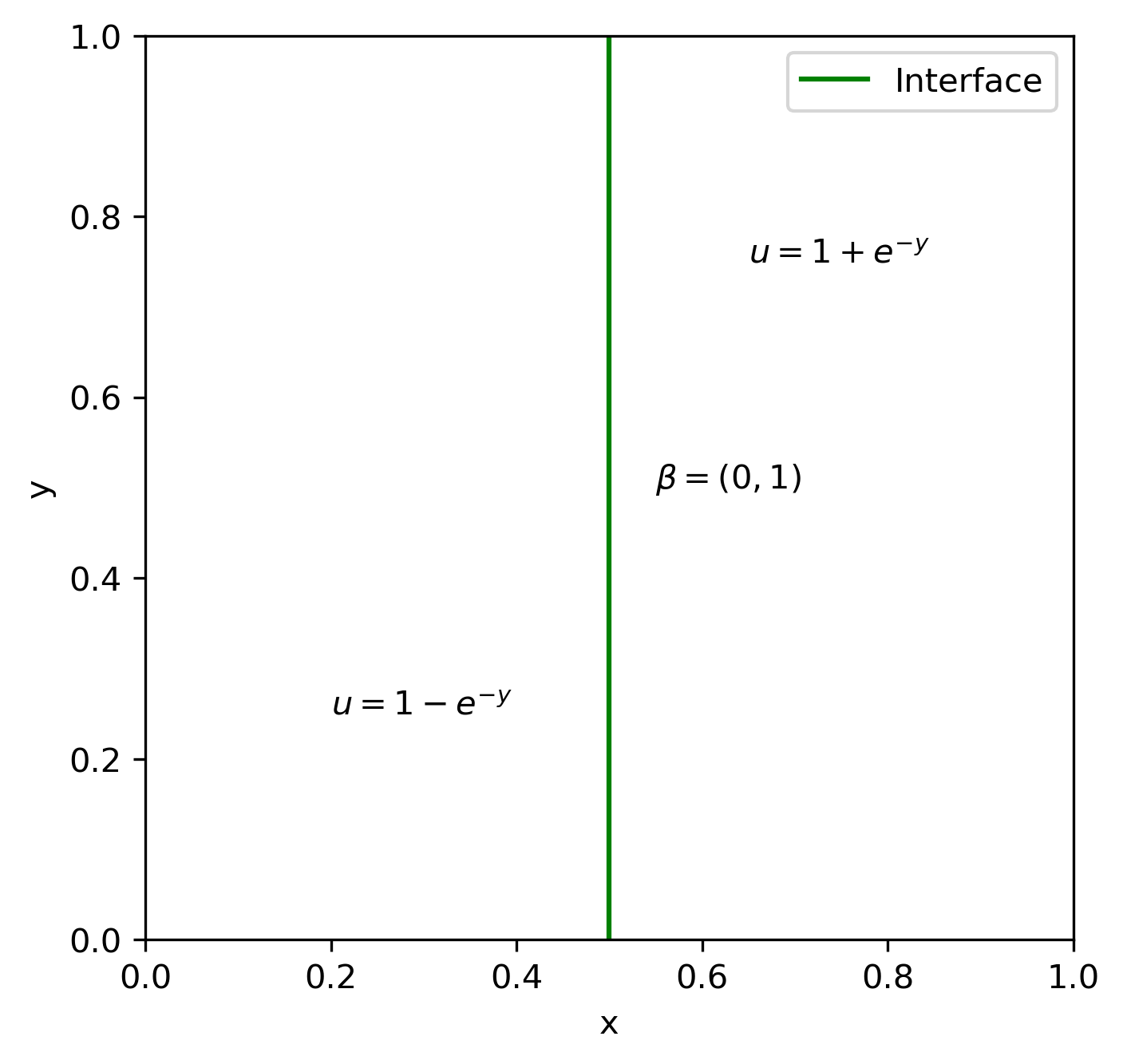}
\end{minipage}%
}%
\hspace{0.2in}
\subfigure[The trace of Figure \ref{comparison2} on $y=0.5$\label{vertical2}]{
\begin{minipage}[t]{0.4\linewidth}
\centering
\includegraphics[width=1.8in]{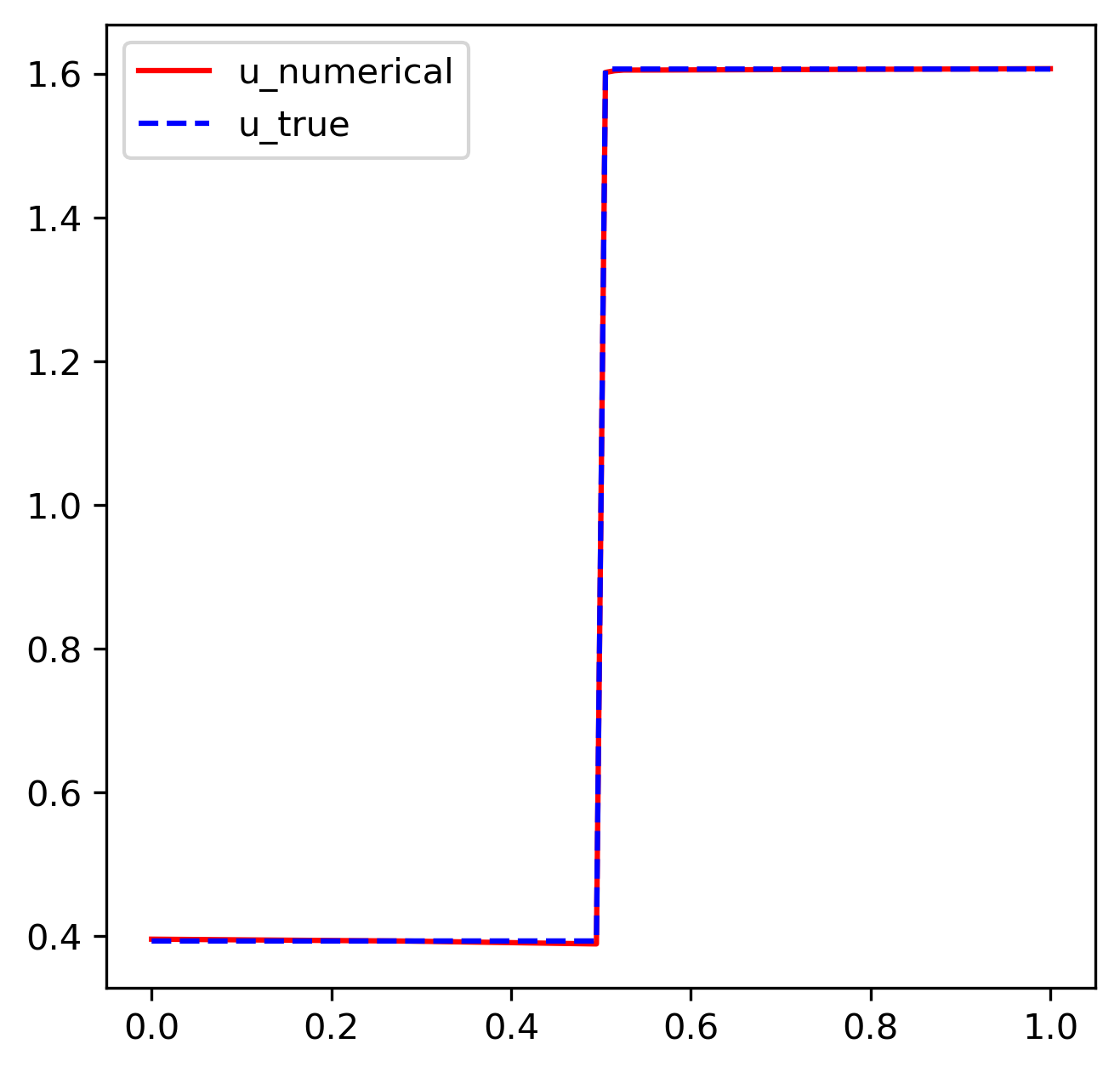}
\end{minipage}%
}%
\\
\subfigure[The exact solution\label{comparison_exact2}]{
\begin{minipage}[t]{0.4\linewidth}
\centering
\includegraphics[width=1.8in]{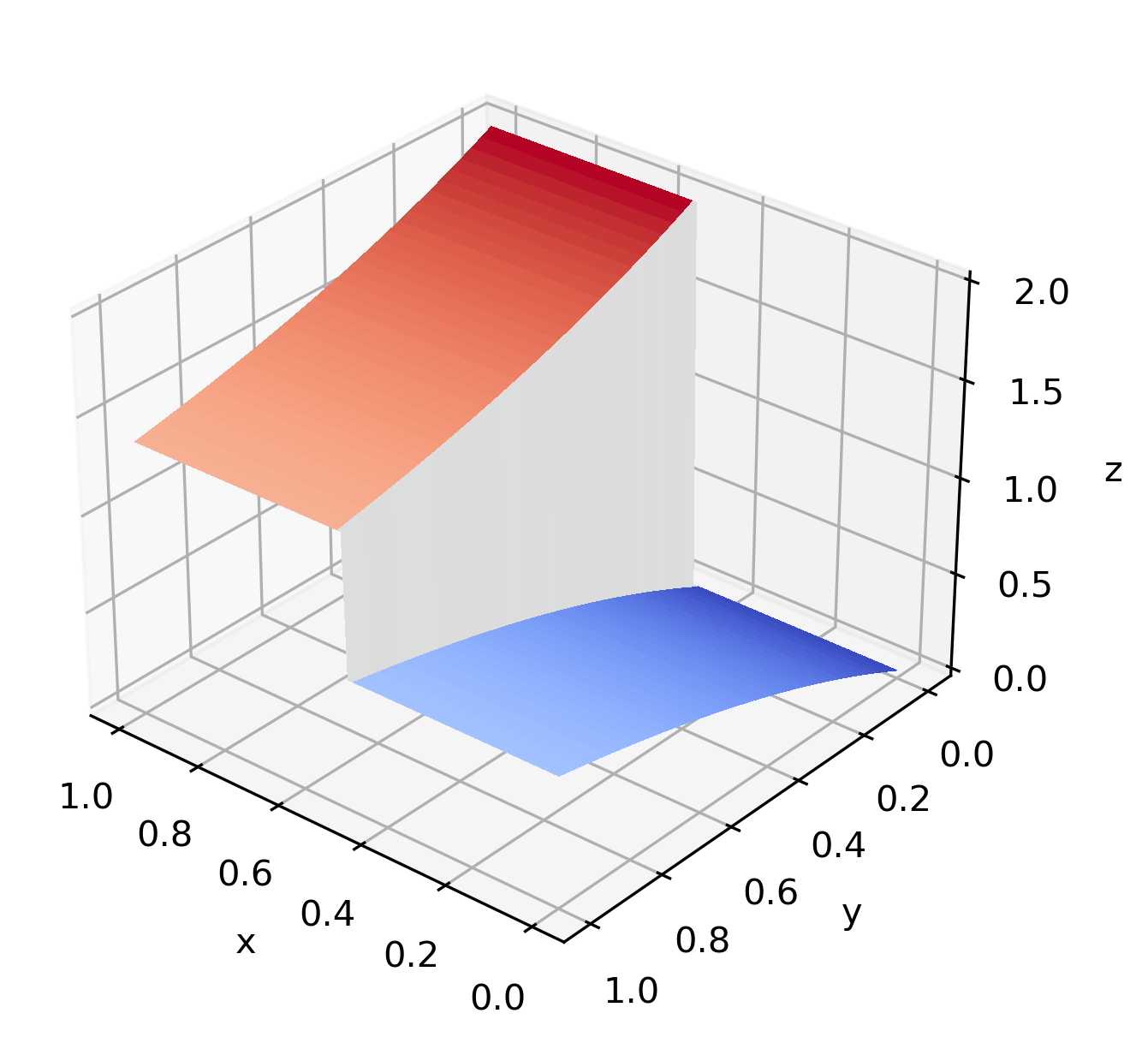}
\end{minipage}%
}%
\hspace{0.2in}
\subfigure[A 2--20--20--1 ReLU NN function approximation\label{comparison2}]{
\begin{minipage}[t]{0.4\linewidth}
\centering
\includegraphics[width=1.8in]{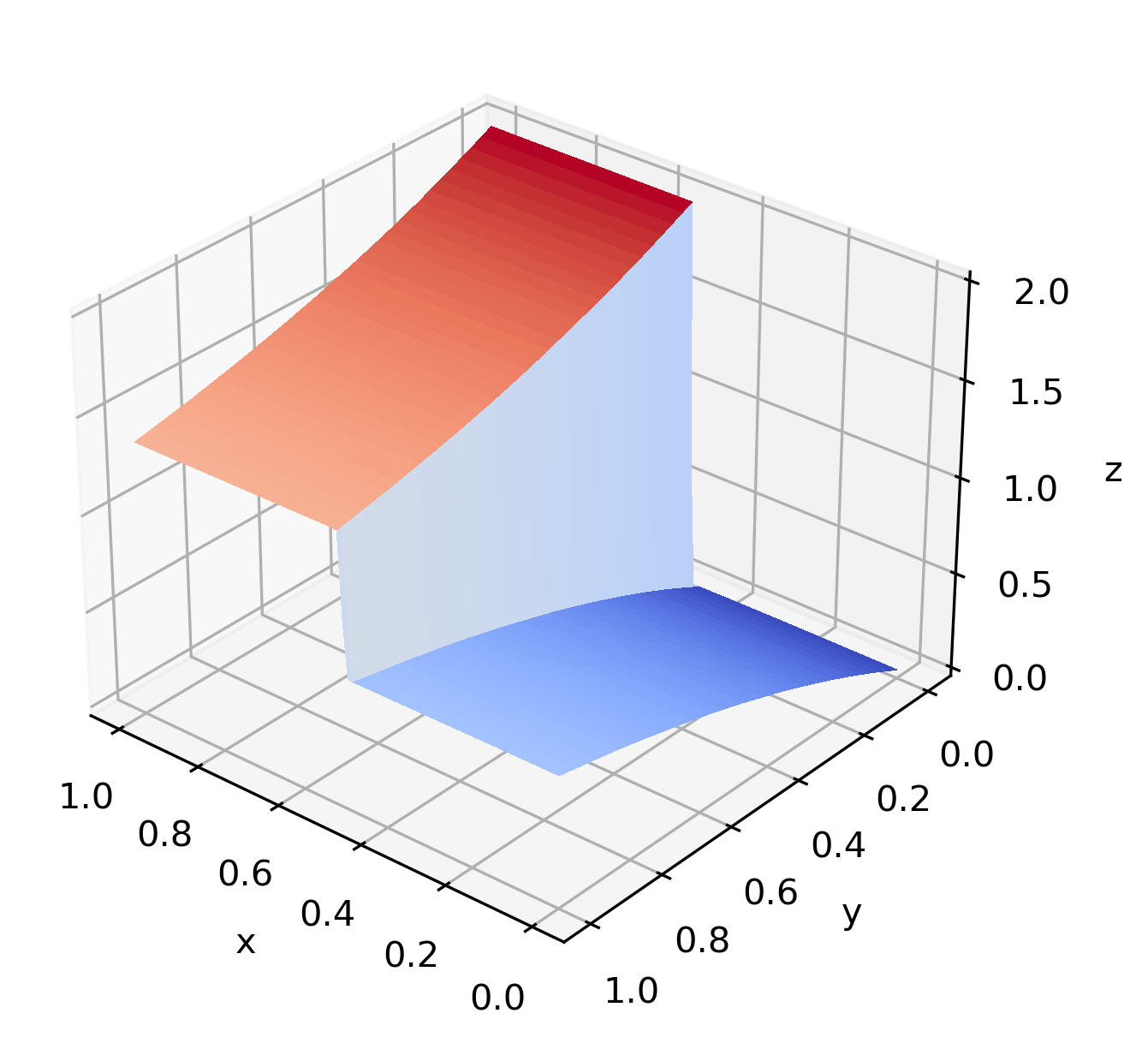}
\end{minipage}%
}%
\\
\subfigure[The breaking hyperplanes of the approximation in Figure \ref{comparison2}\label{breaking2}]{
\begin{minipage}[t]{0.4\linewidth}
\centering
\includegraphics[width=1.8in]{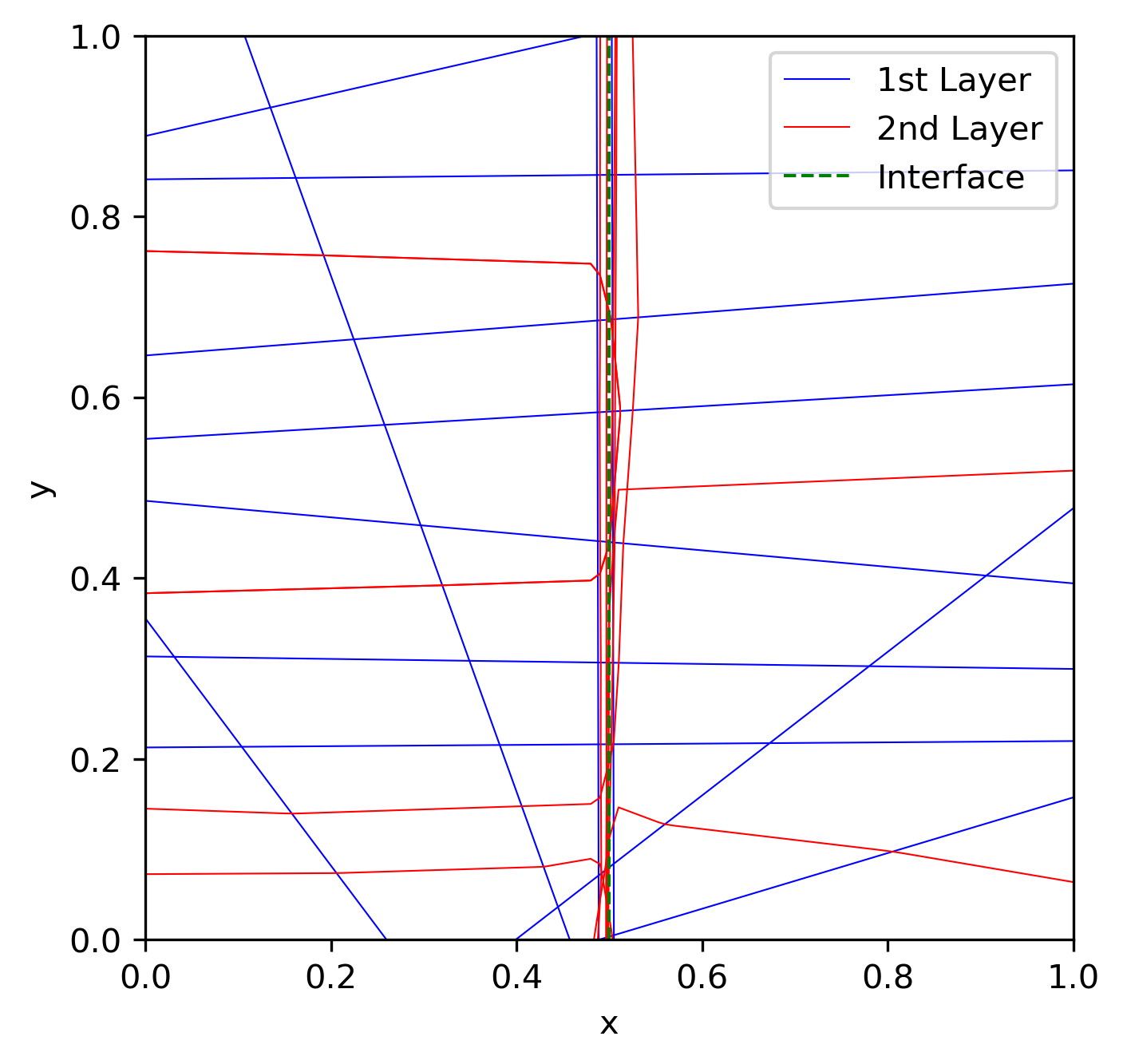}
\end{minipage}%
}%
\caption{Approximation results of the problem in \Cref{test2}}
\end{figure}

The LSNN method was implemented with 50000 iterations for 2--20--20--1 ReLU NN functions. We report the numerical results in \cref{test2 figure,test2 table}. Unlike the previous test problem, the exact solution (\cref{comparison_exact2}) consists of two non-constant smooth parts. The LSNN method is capable of approximating the solution accurately without oscillation (\cref{vertical2,comparison_exact2,comparison2,test2 table}). The 3-layer ReLU NN function approximation has a partition (\cref{breaking2}) similar to that in \Cref{test1} with the second-layer breaking hyperplanes on both sides for approximating the two non-constant smooth parts of the solution.

\begin{table}[htbp]\label{test2 table}
\caption{Relative errors of the problem in \Cref{test2}}
\centering
\begin{tabular}{|l|l|l|l|l|}
\hline
Network structure  &$\frac{\|u-{u}^{_N}_{_\cT}\|_0}{\|u\|_0}$ &$\frac{\vertiii{u-{u}^{_N}_{_\cT}}_{\bm\beta}}{\vertiii{u}_{\bm\beta}}$ & $\frac{\mathcal{L}^{1/2}({u}^{_N}_{_\cT},\bf f)}{\mathcal{L}^{1/2}({u}^{_N}_{_\cT},\bf 0)}$ & Parameters \\ \hline
2--20--20--1  & 0.078036 & 0.013157 & 0.010386   & 501\\ \hline
\end{tabular}
\end{table}

\subsection{A problem with a piecewise smooth inflow boundary}\label{test3}
This example is again a modification of \Cref{test1} by changing the inflow boundary condition to
\begin{equation*}
g(x,y)=\left\{ \begin{array}{rl}
 1-\sin(2\pi x),& (x,y)\in \Gamma^1_-\equiv \{(x,0): x\in(0,1/2)\}, \\[2mm]
 5/2-x, &(x,y)\in \Gamma^2_-=\Gamma_-\setminus \Gamma_-^1.
\end{array}\right.
\end{equation*} 
The exact solution of this test problem is
\begin{equation}
u(x,y)=\left\{ \begin{array}{rl}
 1-\sin(2\pi x)e^{-y},& (x,y)\in \Omega_1, \\[2mm]
 1+(3/2-x)e^{-y}, & (x,y)\in \Omega_2.
\end{array}\right.
\end{equation}

\begin{figure}[htbp]\label{test3 figure}
\centering
\subfigure[The interface\label{interface3}]{
\begin{minipage}[t]{0.4\linewidth}
\centering
\includegraphics[width=1.8in]{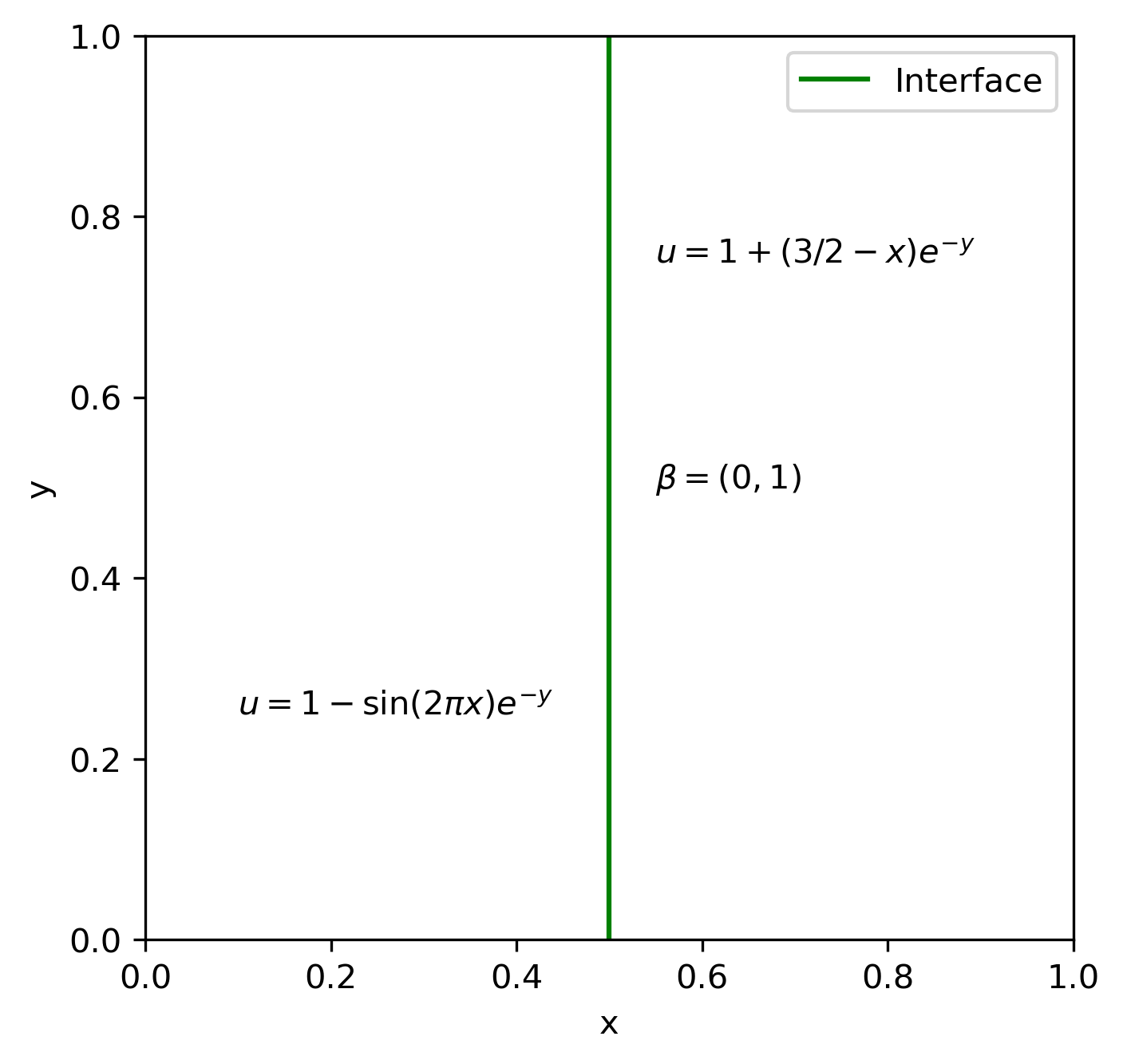}
\end{minipage}%
}%
\hspace{0.2in}
\subfigure[The trace of Figure \ref{comparison3} on $y=0.5$\label{vertical32}]{
\begin{minipage}[t]{0.4\linewidth}
\centering
\includegraphics[width=1.8in]{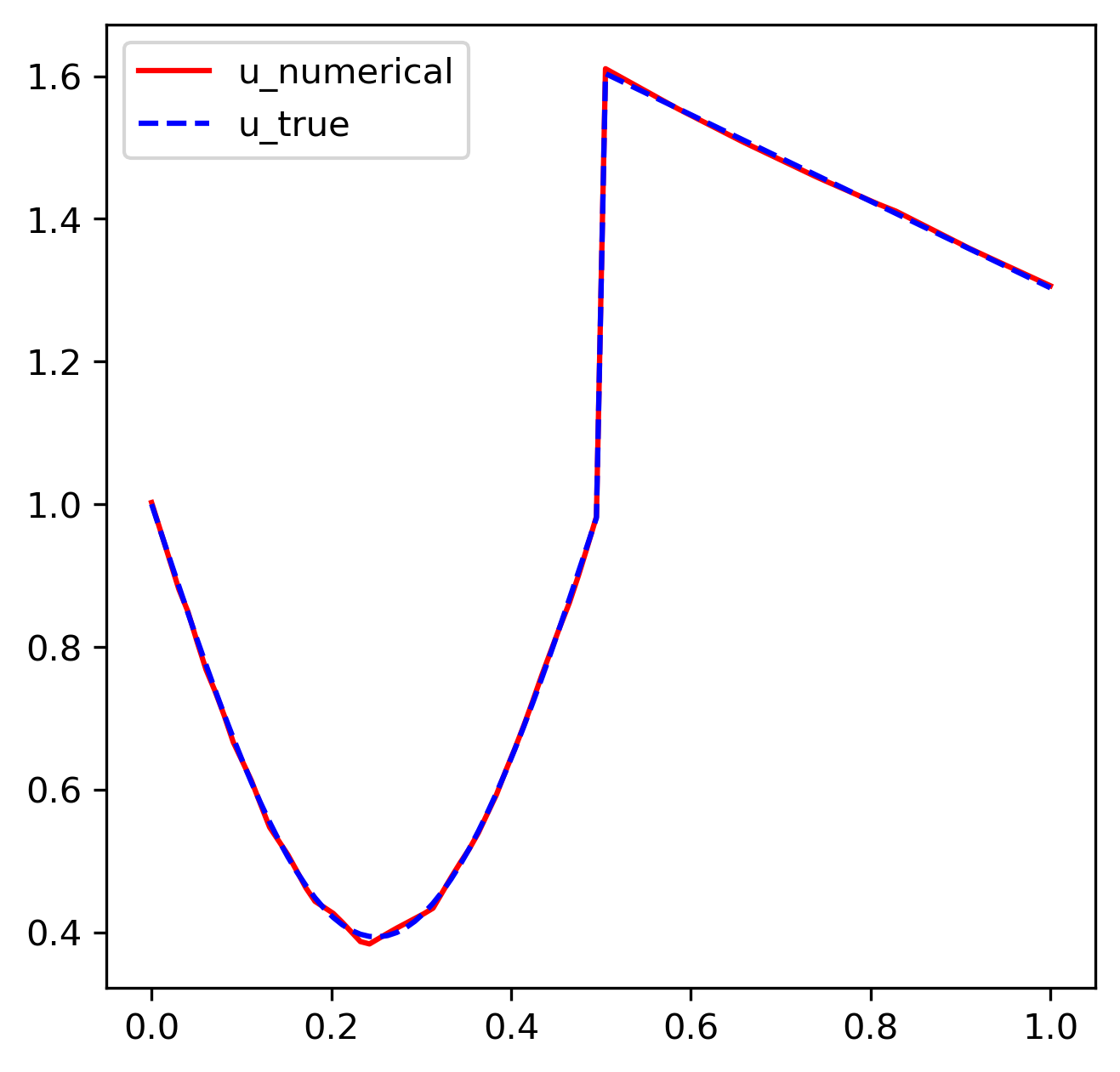}
\end{minipage}%
}%
\\
\subfigure[The exact solution\label{comparison_exact3}]{
\begin{minipage}[t]{0.4\linewidth}
\centering
\includegraphics[width=1.8in]{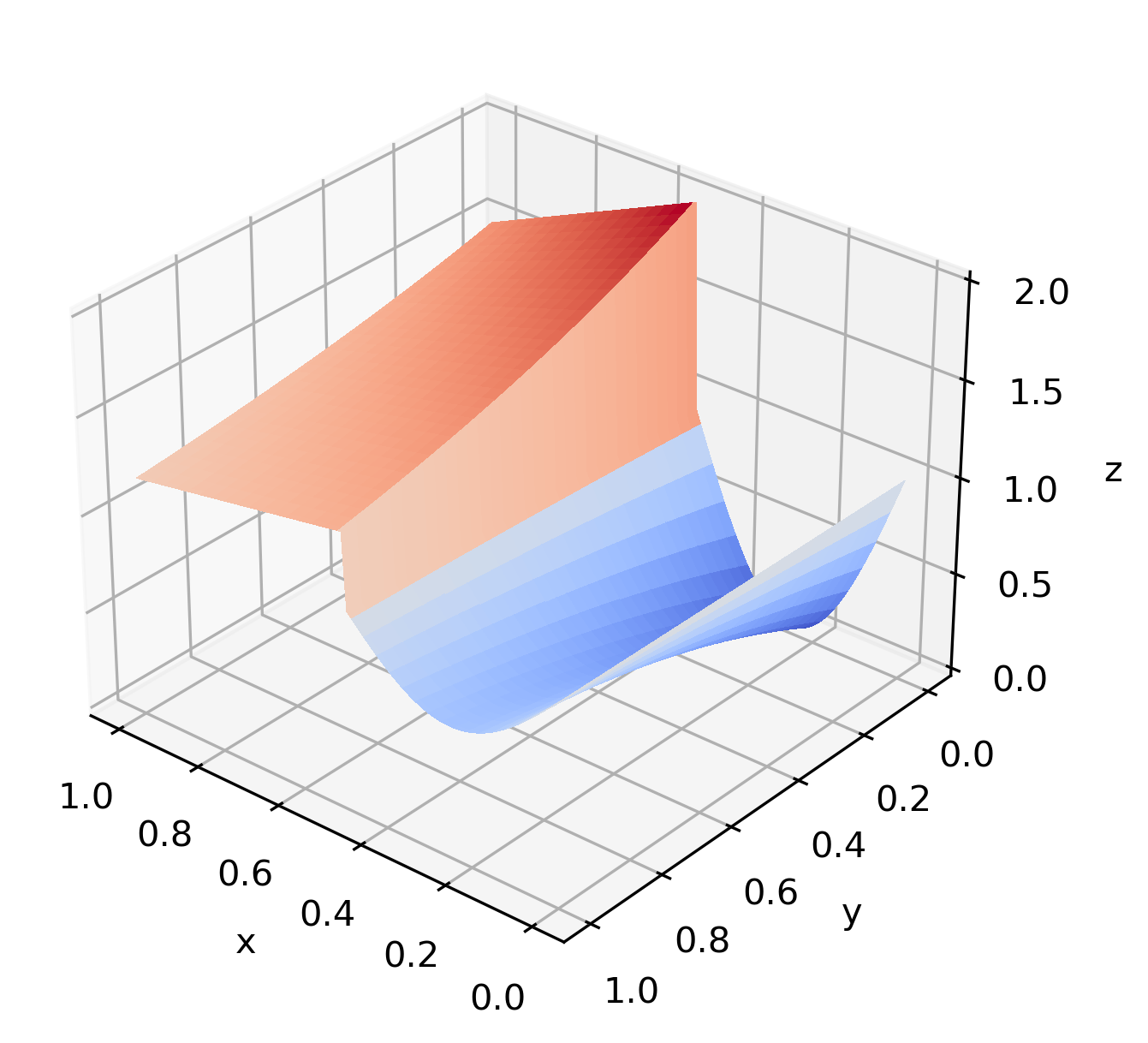}
\end{minipage}%
}%
\hspace{0.2in}
\subfigure[A 2--40--40--1 ReLU NN function approximation\label{comparison3}]{
\begin{minipage}[t]{0.4\linewidth}
\centering
\includegraphics[width=1.8in]{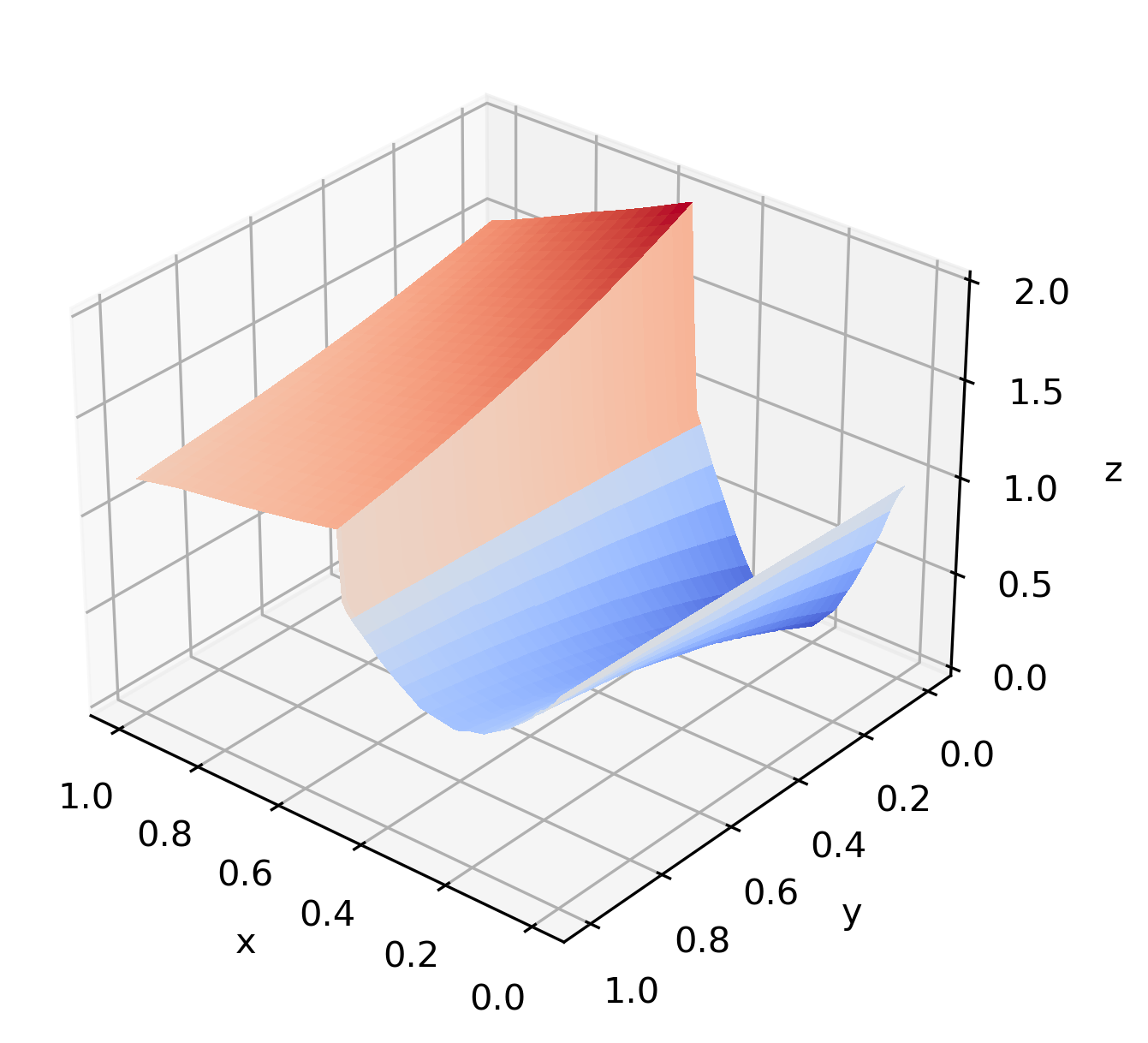}
\end{minipage}%
}%
\\
\subfigure[The breaking hyperplanes of the approximation in Figure \ref{comparison3}\label{breaking32}]{
\begin{minipage}[t]{0.4\linewidth}
\centering
\includegraphics[width=1.8in]{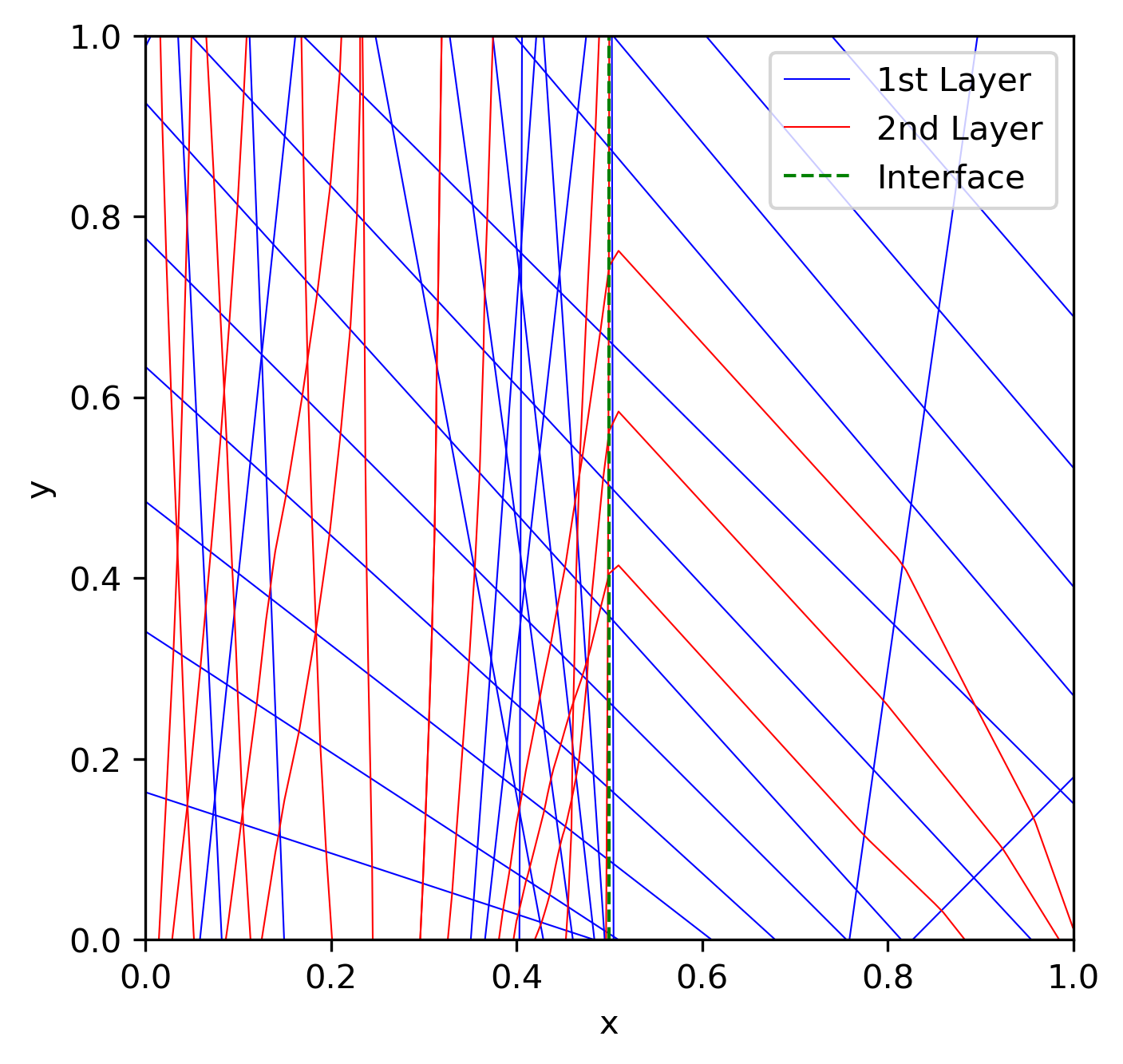}
\end{minipage}%
}%
\caption{Approximation results of the problem in \Cref{test3}}
\end{figure}

The LSNN method was implemented with 100000 iterations for 2--40--40--1 ReLU NN functions. We report the numerical results in \cref{test3 figure,test3 table}. Since the solution on the inflow boundary consists of two non-constant smooth curves, we increased the size of the neural network to obtain a more accurate solution. \cref{comparison_exact3,comparison3,test3 table} show that the approximation is accurate pointwisely and in average. The traces (\cref{vertical32}) on $y=0.5$ exhibit no oscillation and we note a few corners on the curve, verifying that the ReLU NN function approximation is a CPWL function. The partition generated by the breaking hyperplanes (\cref{breaking32}) of the approximation shows how the exact solution was approximated.

\begin{table}[htbp]\label{test3 table}
\caption{Relative errors of the problem in \Cref{test3}}
\centering
\begin{tabular}{|l|l|l|l|l|}
\hline
Network structure  &$\frac{\|u-{u}^{_N}_{_\cT}\|_0}{\|u\|_0}$ &$\frac{\vertiii{u-{u}^{_N}_{_\cT}}_{\bm\beta}}{\vertiii{u}_{\bm\beta}}$ & $\frac{\mathcal{L}^{1/2}({u}^{_N}_{_\cT},\bf f)}{\mathcal{L}^{1/2}({u}^{_N}_{_\cT},\bf 0)}$ & Parameters \\ \hline
2--40--40--1  & 0.041491 & 0.016480 & 0.012733   & 1801\\ \hline
\end{tabular}
\end{table}

\subsection{A problem with a piecewise constant advection velocity field}\label{test4}
Let $\bar{\Omega}=\bar{\Upsilon}_1\cup \bar{\Upsilon}_2$ and 
\[
    \Upsilon_1=\{(x,y)\in\Omega:y\ge x\}\text{ and }
    \Upsilon_2=\Omega\setminus\Upsilon_1.
\]

This test problem has a piecewise constant advective velocity field given by
\begin{equation}
\bm{\beta}(x,y) =\left\{ \begin{array}{rl}
(-1,\sqrt{2}-1)^T,&(x,y)\in\Upsilon_1,\\[2mm]
(1-\sqrt{2},1)^T,&(x,y)\in\Upsilon_2,
\end{array}\right.
\end{equation}
and the inflow boundary $\Gamma_{-}=\{(1,y):y\in(0,1)\}\cup\{(x,0):x\in(0,1)\}$. For the inflow boundary condition,
\begin{equation*}
g(x,y)=\left\{ \begin{array}{rl}
 x\exp(x/(\sqrt{2}-1)),& (x,y)\in \Gamma^1_-\equiv \{(x,0): x\in(0,1)\}, \\[2mm]
 (11+(\sqrt{2}-1)y)\exp(1/(\sqrt{2}-1)), &(x,y)\in \Gamma^2_-=\Gamma_-\setminus \Gamma_-^1,
\end{array}\right.
\end{equation*} 

the exact solution of this test problem is
\begin{equation}
u(x,y) =\left\{ \begin{array}{rl}
(y+(\sqrt{2}-1)x)\exp(\sqrt{2}x+y),&(x,y)\in\widehat{\Upsilon}_{11},\\[2mm]
(y+(\sqrt{2}-1)x+10)\exp(\sqrt{2}x+y),&(x,y)\in\widehat{\Upsilon}_{12},\\[2mm]
(x+(\sqrt{2}-1)y)\exp(x/(\sqrt{2}-1)),&(x,y)\in\widehat{\Upsilon}_{21},\\[2mm]
(x+(\sqrt{2}-1)y+10)\exp(x/(\sqrt{2}-1)),&(x,y)\in\widehat{\Upsilon}_{22},
\end{array}\right.
\end{equation}

where
\begin{multline*}
\widehat{\Upsilon}_{11}=\{(x,y)\in\Upsilon_1: y<(1-\sqrt{2})x+1\},\,\, \widehat{\Upsilon}_{12}=\Upsilon_1\setminus\widehat{\Upsilon}_{11},\\[2mm]
\widehat{\Upsilon}_{21}=\{(x,y)\in\Upsilon_2:y<\tfrac{1}{1-\sqrt{2}}(x-1)\},\text{ and }\widehat{\Upsilon}_{22}=\Upsilon_2\setminus\widehat{\Upsilon}_{21}.
\end{multline*}

\begin{figure}[htbp]\label{test4 figure}
\centering
\subfigure[The interface\label{interface42}]{
\begin{minipage}[t]{0.4\linewidth}
\centering
\includegraphics[width=1.8in]{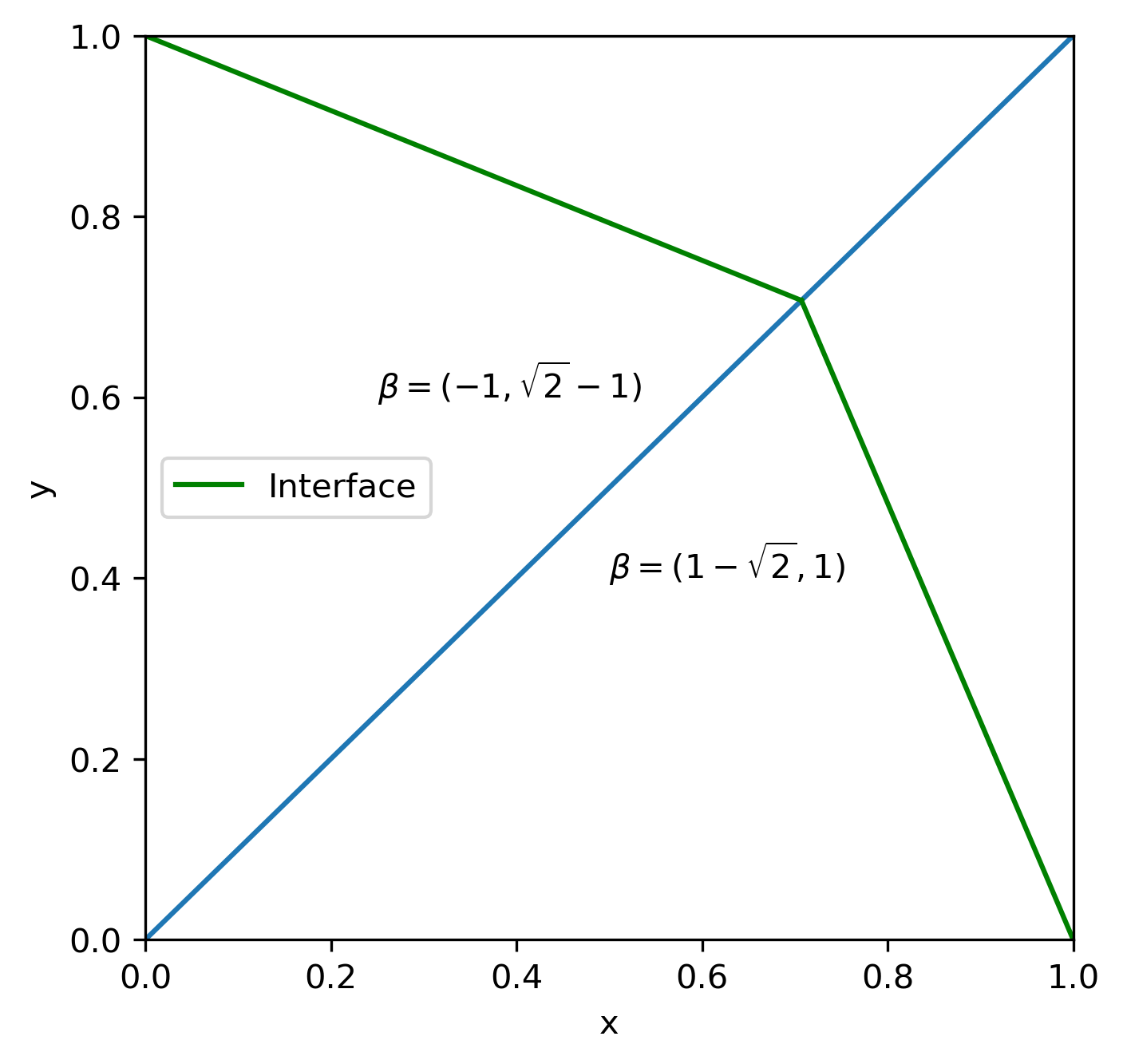}
\end{minipage}%
}%
\hspace{0.2in}
\subfigure[The exact solution\label{comparison_exact4}]{
\begin{minipage}[t]{0.4\linewidth}
\centering
\includegraphics[width=1.8in]{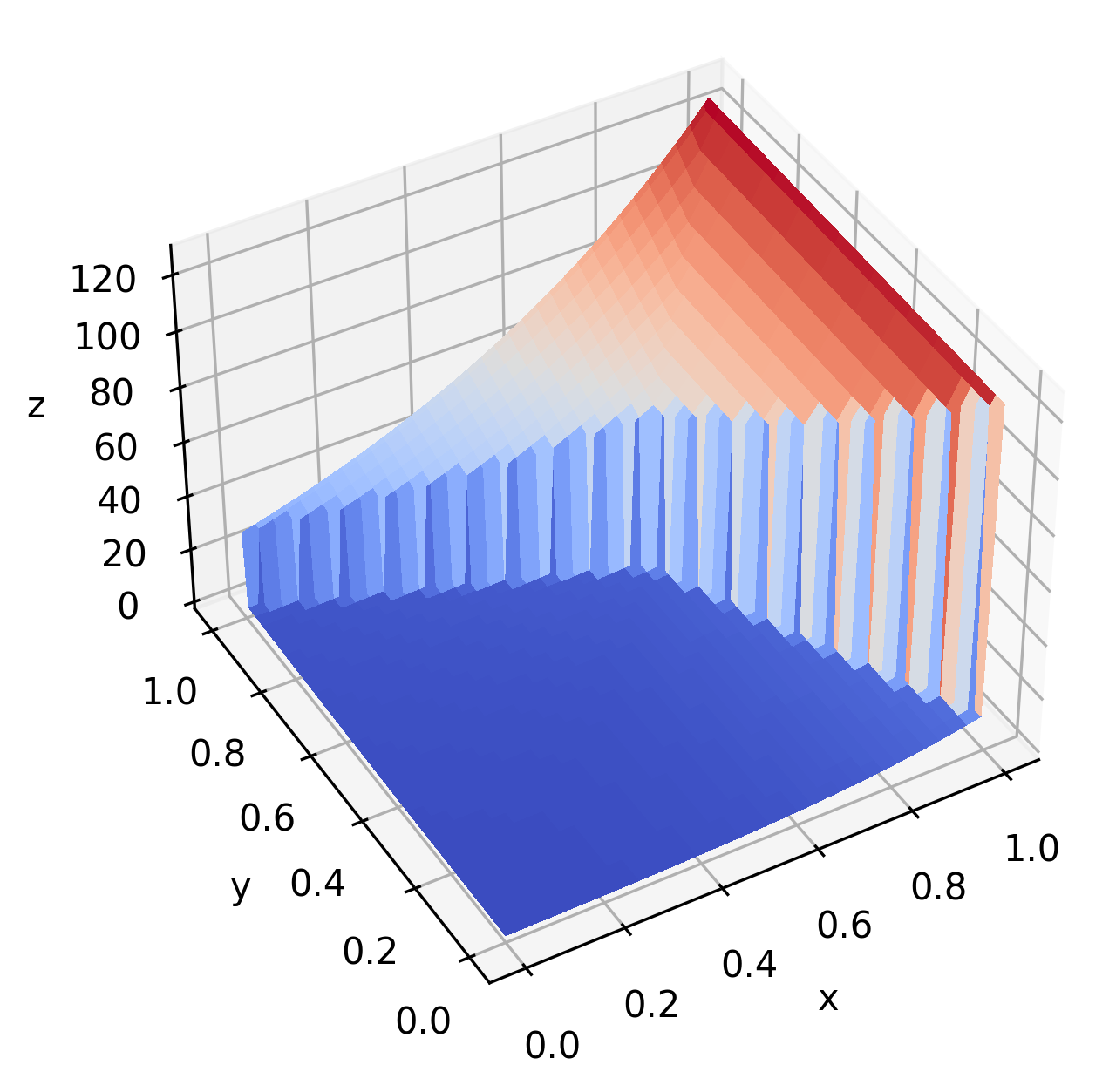}
\end{minipage}%
}%
\\
\subfigure[A 2--450--1 ReLU NN function approximation\label{comparison4-2}]{
\begin{minipage}[t]{0.4\linewidth}
\centering
\includegraphics[width=1.8in]{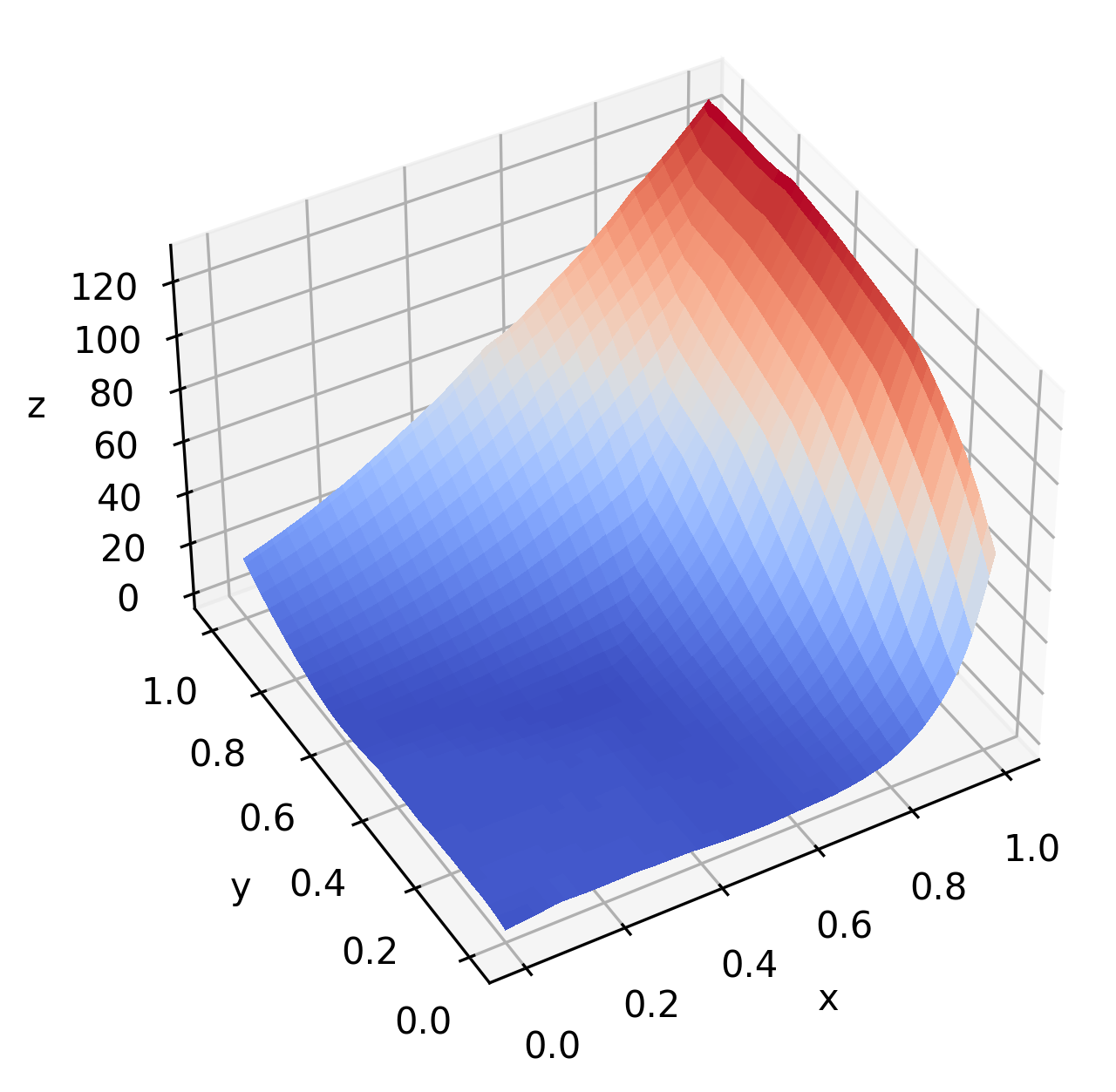}
\end{minipage}%
}%
\hspace{0.2in}
\subfigure[A 2--40--40--1 ReLU NN function approximation\label{comparison4}]{
\begin{minipage}[t]{0.4\linewidth}
\centering
\includegraphics[width=1.8in]{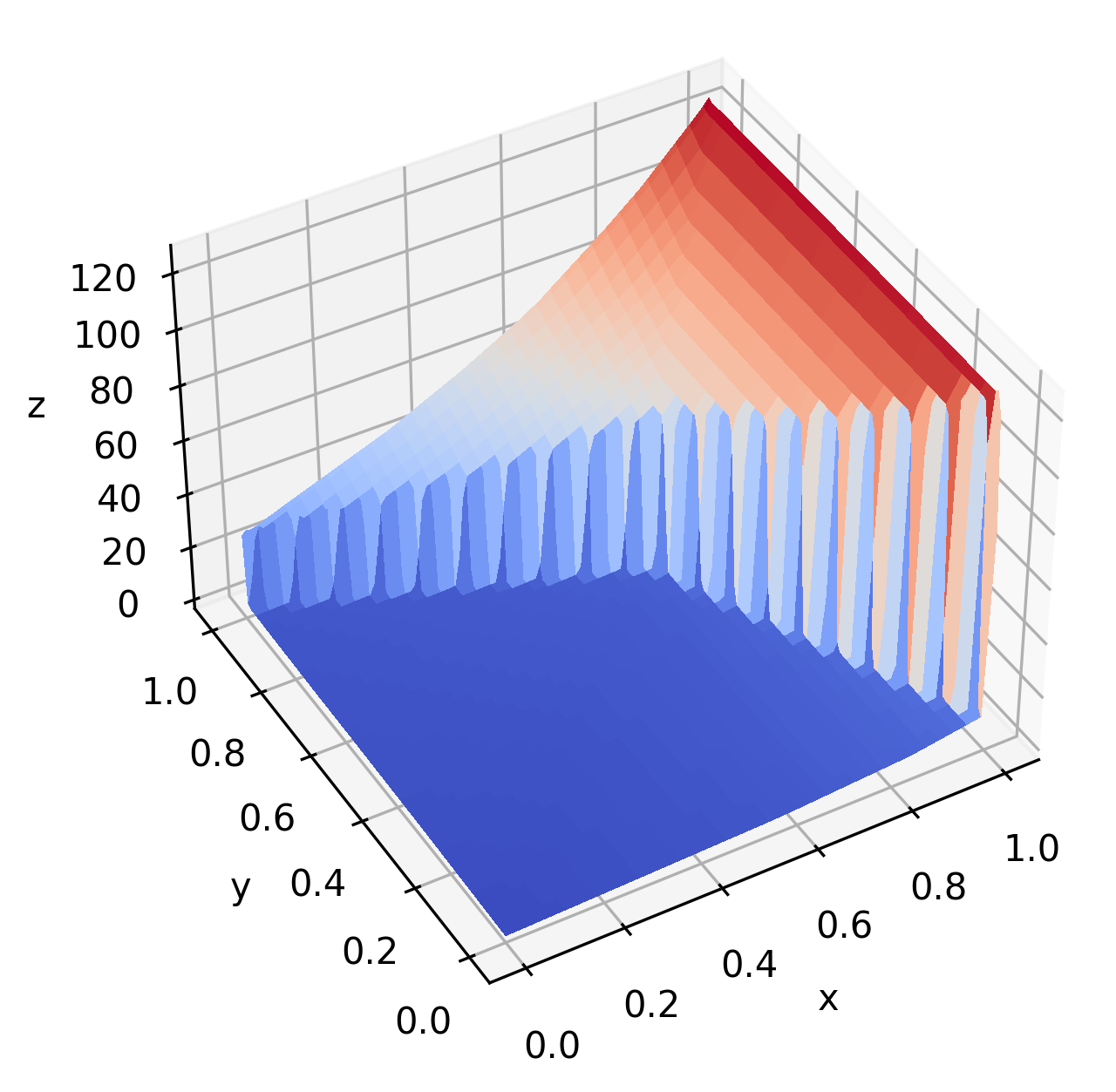}
\end{minipage}%
}%
\\
\subfigure[The trace of Figure \ref{comparison4-2} on $y=x$\label{vertical4-2}]{
\begin{minipage}[t]{0.4\linewidth}
\centering
\includegraphics[width=1.8in]{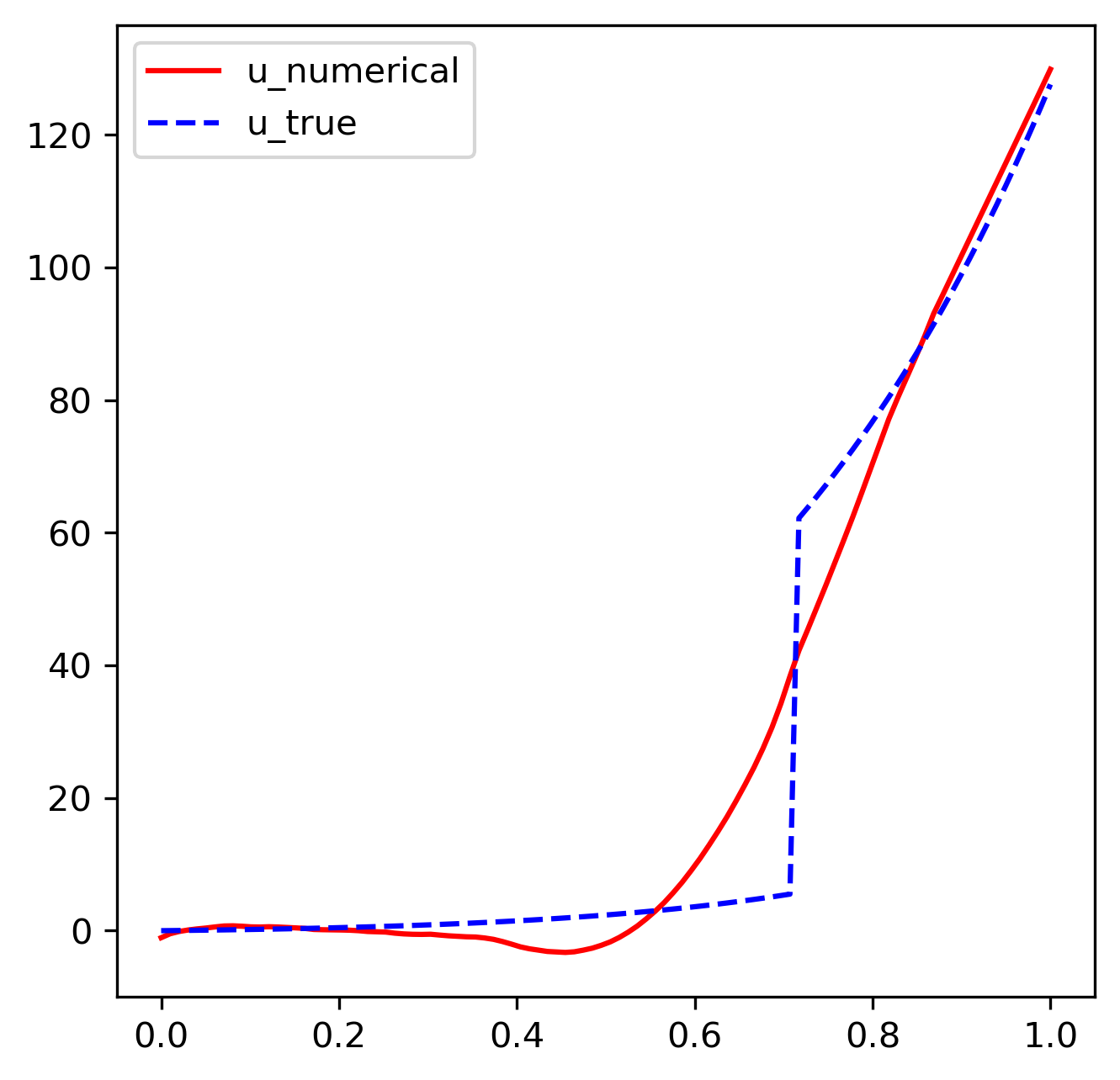}
\end{minipage}%
}%
\hspace{0.2in}
\subfigure[The trace of Figure \ref{comparison4} on $y=x$\label{vertical42}]{
\begin{minipage}[t]{0.4\linewidth}
\centering
\includegraphics[width=1.8in]{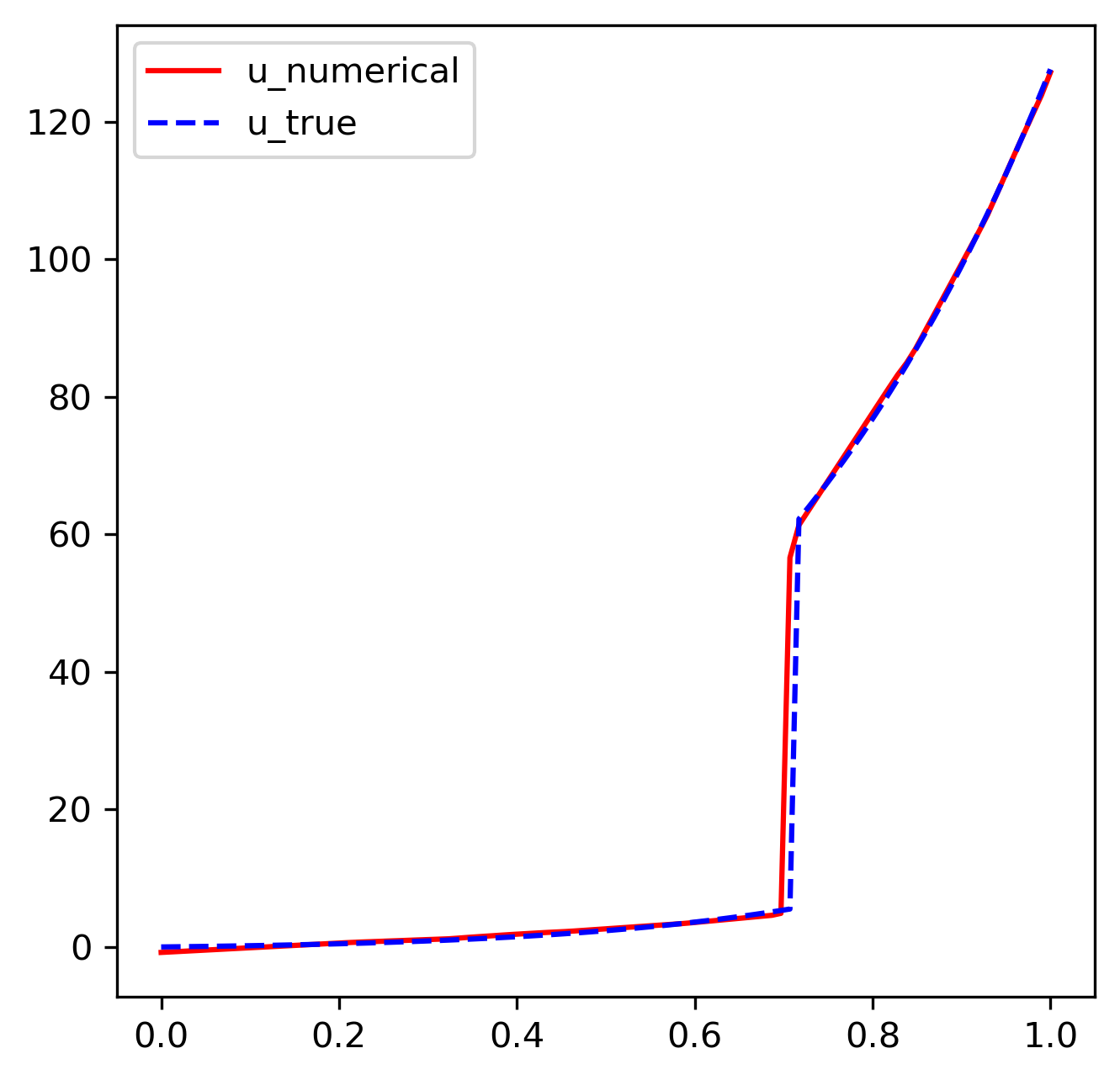}
\end{minipage}%
}%
\\
\subfigure[The breaking hyperplanes of the approximation in Figure \ref{comparison4-2}\label{breaking4-2}]{
\begin{minipage}[t]{0.4\linewidth}
\centering
\includegraphics[width=1.8in]{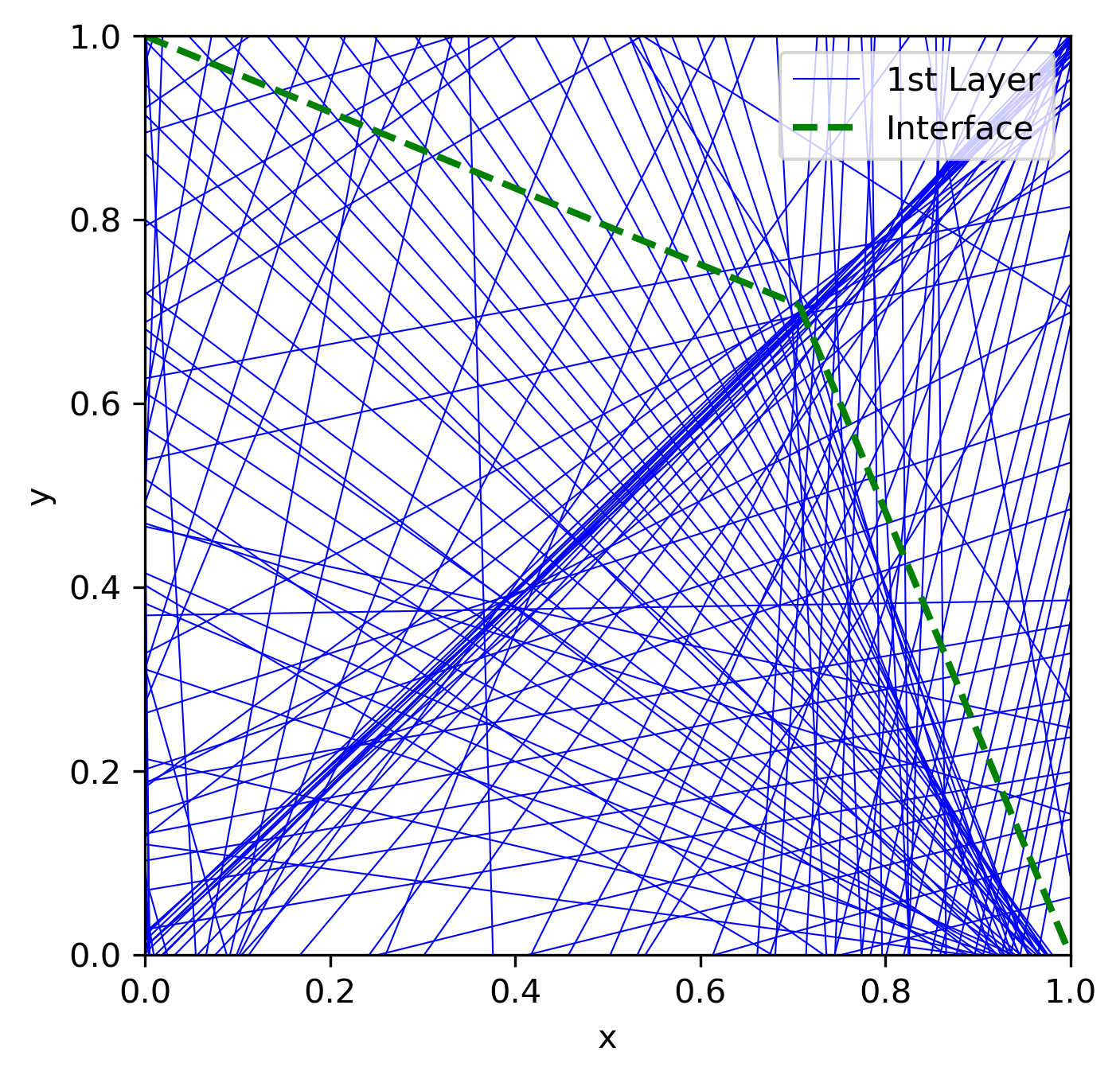}
\end{minipage}%
}%
\hspace{0.2in}
\subfigure[The breaking hyperplanes of the approximation in Figure \ref{comparison4}\label{breaking42}]{
\begin{minipage}[t]{0.4\linewidth}
\centering
\includegraphics[width=1.8in]{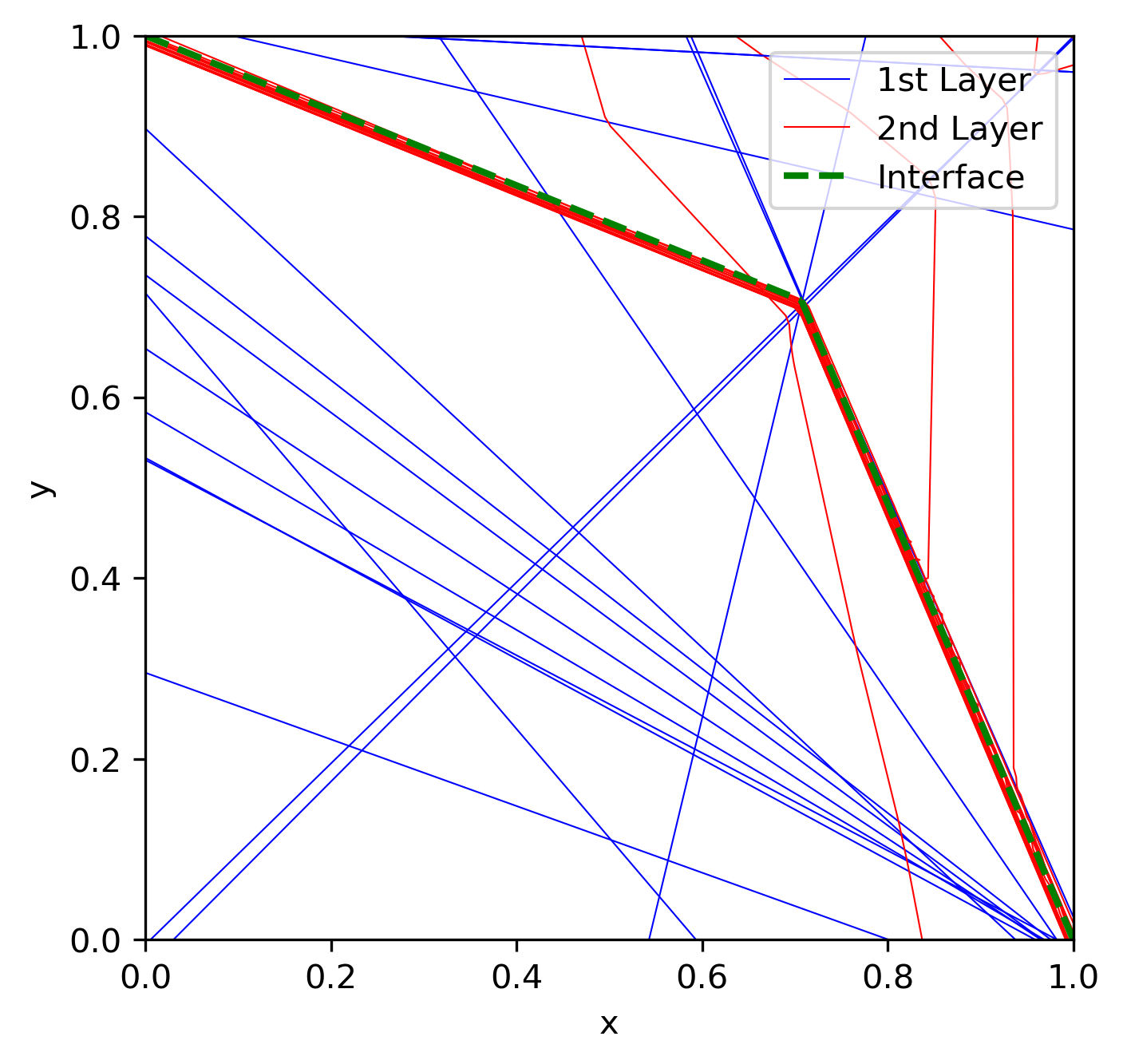}
\end{minipage}%
}%
\caption{Approximation results of the problem in \Cref{test4}}
\end{figure}

The LSNN method was implemented with 300000 iterations for 2--450--1 and 2--40--40--1 ReLU NN functions. We report the numerical results in \cref{test4 figure,test4 table}. This example compares the approximation differences with the one-hidden-layer NN (a universal approximator) with the same number of degrees of freedom. The exact solution (\cref{comparison_exact4}) consists of four non-constant smooth parts and has a non-constant jump along two connected line segments (\cref{interface42}). The traces (\cref{vertical42}) of the exact and numerical solutions, the approximation (\cref{comparison4}), and \cref{test4 table} indicate that the 3-layer ReLU NN function approximation is accurate pointwisely and in average. Most of the second-layer breaking hyperplanes (\cref{breaking42}) are along the discontinuous interface, which correspond to the sharp transition layer of the approximation for the jump. On the other hand, \cref{comparison4-2,vertical4-2,breaking4-2,test4 table} show that the one-hidden-layer NN failed to approximate the solution, especially around the interface.

\begin{table}[htbp]\label{test4 table}
\caption{Relative errors of the problem in \Cref{test4}}
\centering
\begin{tabular}{|l|l|l|l|l|}
\hline
Network structure  &$\frac{\|u-{u}^{_N}_{_\cT}\|_0}{\|u\|_0}$ &$\frac{\vertiii{u-{u}^{_N}_{_\cT}}_{\bm\beta}}{\vertiii{u}_{\bm\beta}}$ & $\frac{\mathcal{L}^{1/2}({u}^{_N}_{_\cT},\bf f)}{\mathcal{L}^{1/2}({u}^{_N}_{_\cT},\bf 0)}$ & Parameters \\ \hline
2--40--40--1  & 0.058884 & 0.073169 & 0.038245   & 1801\\ \hline
2--450--1  & 0.243956 & 0.258038 & 0.225756   & 1801\\ \hline
\end{tabular}
\end{table}

\subsection{A problem with a variable advection velocity field}\label{test5}
This test problem has the variable advective velocity field $\bm{\beta}(x,y) = (1,2x),\,\, (x,y)\in\Omega$, and the inflow boundary of the problem is $\Gamma_{-}=\{(0,y):y\in(0,1)\}\cup\{(x,0):x\in(0,1)\}$. The inflow boundary condition is given by
\begin{equation*}
g(x,y)=\left\{ \begin{array}{rl}
 y+2,& (x,y)\in \Gamma^1_-\equiv\{(0,y): y\in[\frac{1}{5},1)\}, \\[2mm]
 (y-x^2)e^{-x}, & (x,y)\in \Gamma^2_-=\Gamma_-\setminus \Gamma_-^1.
\end{array}\right.
\end{equation*} 
The exact solution of this test problem is
\begin{equation}
u(x,y)=\left\{ \begin{array}{rl}
 (y-x^2)e^{-x},& (x,y)\in \Omega_1\equiv\{(x,y)\in\Omega:y< x^2+\frac{1}{5}\}, \\[2mm]
 (y-x^2+2)e^{-x}, & (x,y)\in \Omega_2=\Omega\setminus\Omega_1.
\end{array}\right.
\end{equation}

\begin{figure}[htbp]\label{test5 figure}
\centering
\subfigure[The interface\label{interface5}]{
\begin{minipage}[t]{0.4\linewidth}
\centering
\includegraphics[width=1.8in]{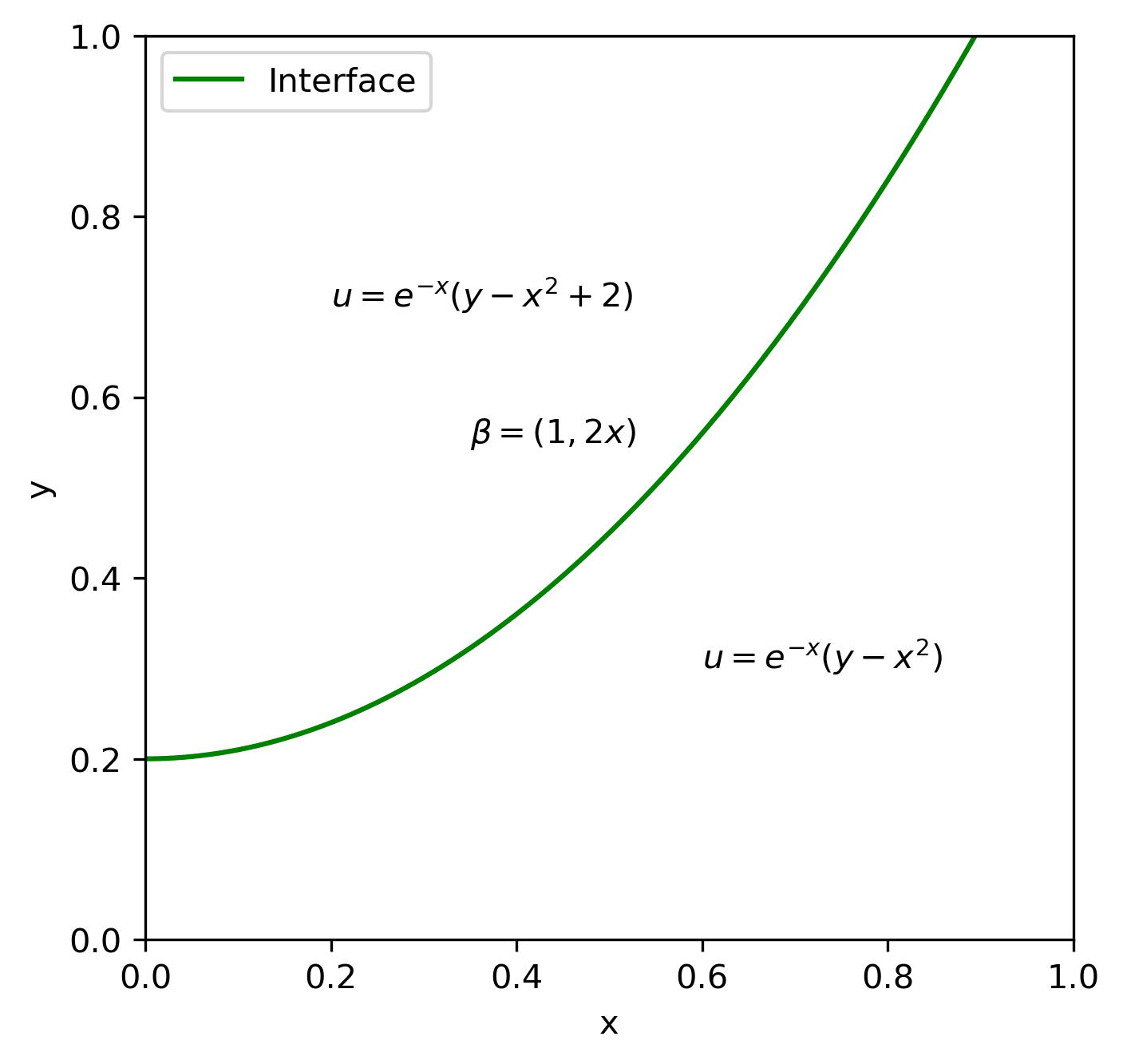}
\end{minipage}%
}%
\hspace{0.2in}
\subfigure[The trace of Figure \ref{comparison5} on $y=1-x$\label{vertical5}]{
\begin{minipage}[t]{0.4\linewidth}
\centering
\includegraphics[width=1.8in]{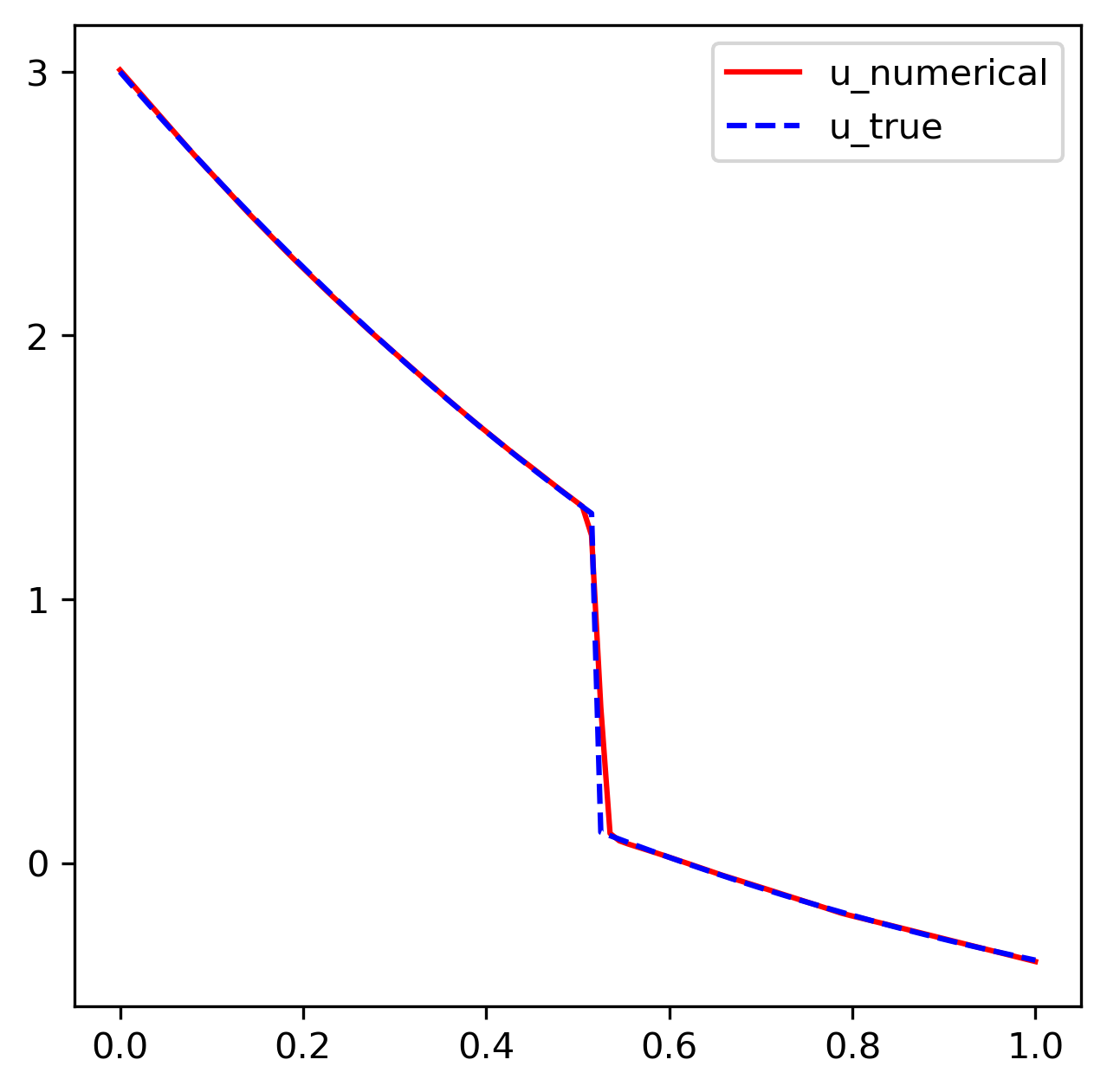}
\end{minipage}%
}%
\\
\subfigure[The exact solution\label{comparison_exact5}]{
\begin{minipage}[t]{0.4\linewidth}
\centering
\includegraphics[width=1.8in]{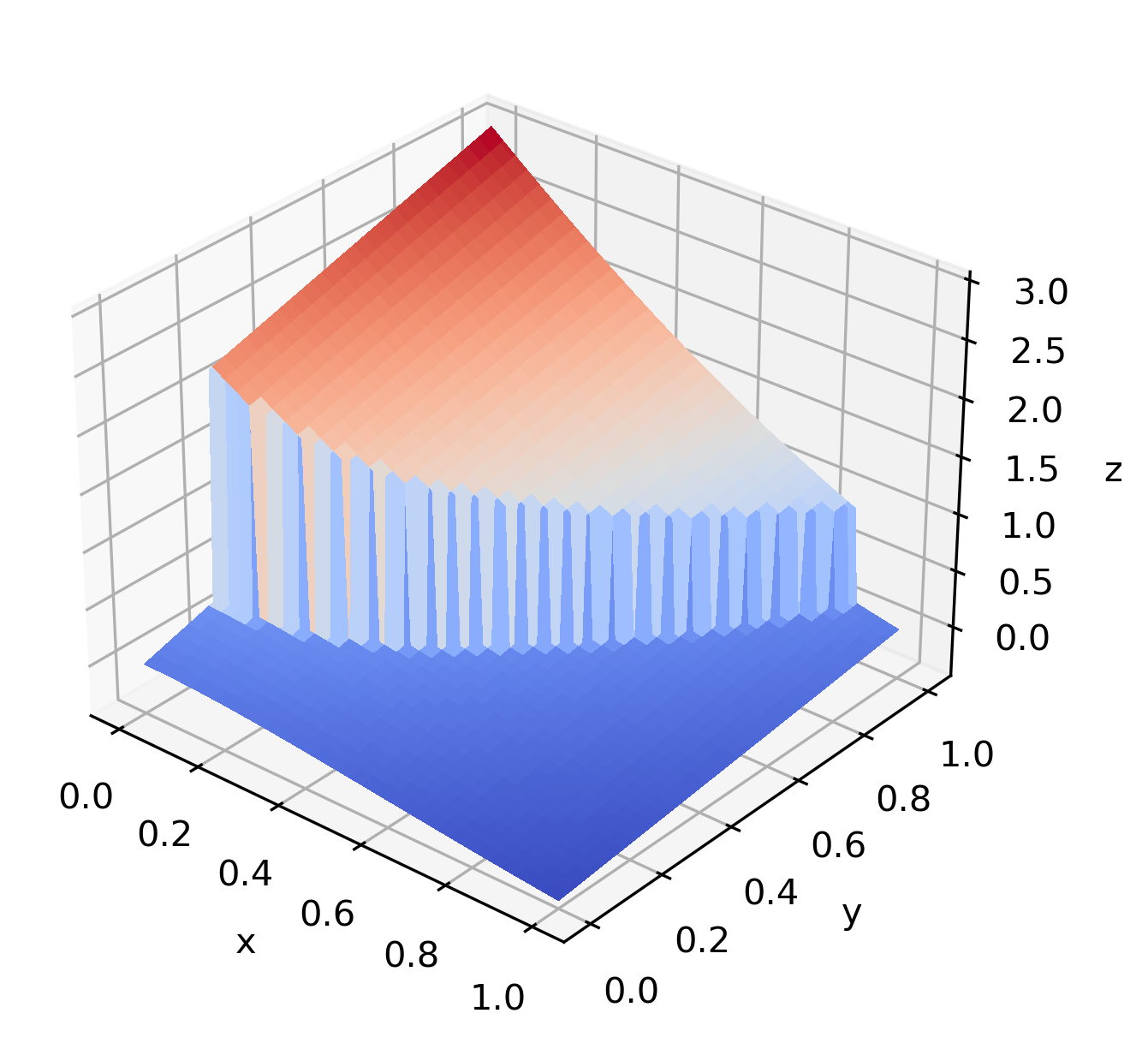}
\end{minipage}%
}%
\hspace{0.2in}
\subfigure[A 2--60--60--1 ReLU NN function approximation\label{comparison5}]{
\begin{minipage}[t]{0.4\linewidth}
\centering
\includegraphics[width=1.8in]{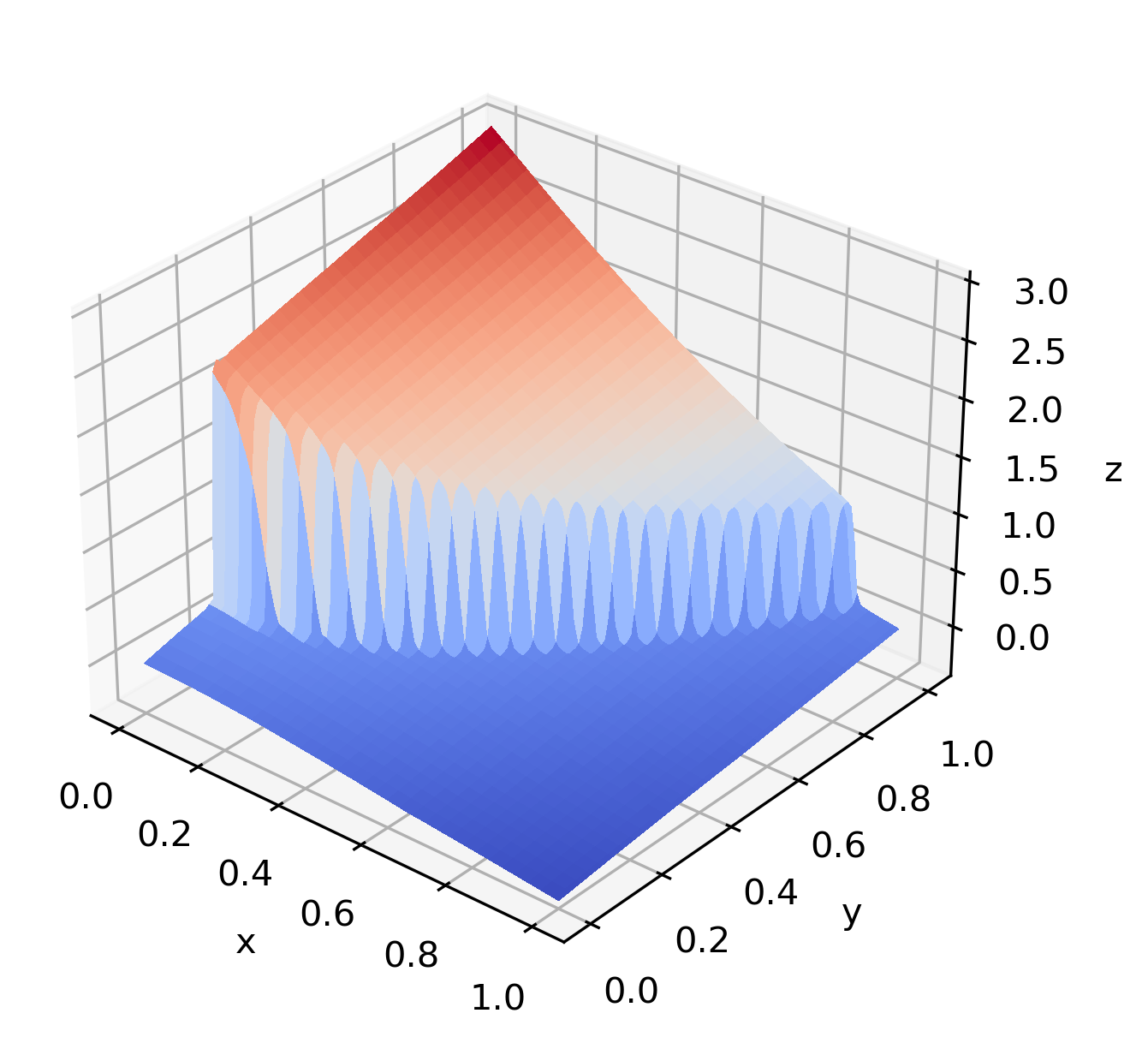}
\end{minipage}%
}%
\\
\subfigure[The breaking hyperplanes of the approximation in Figure \ref{comparison5}\label{breaking5}]{
\begin{minipage}[t]{0.4\linewidth}
\centering
\includegraphics[width=1.8in]{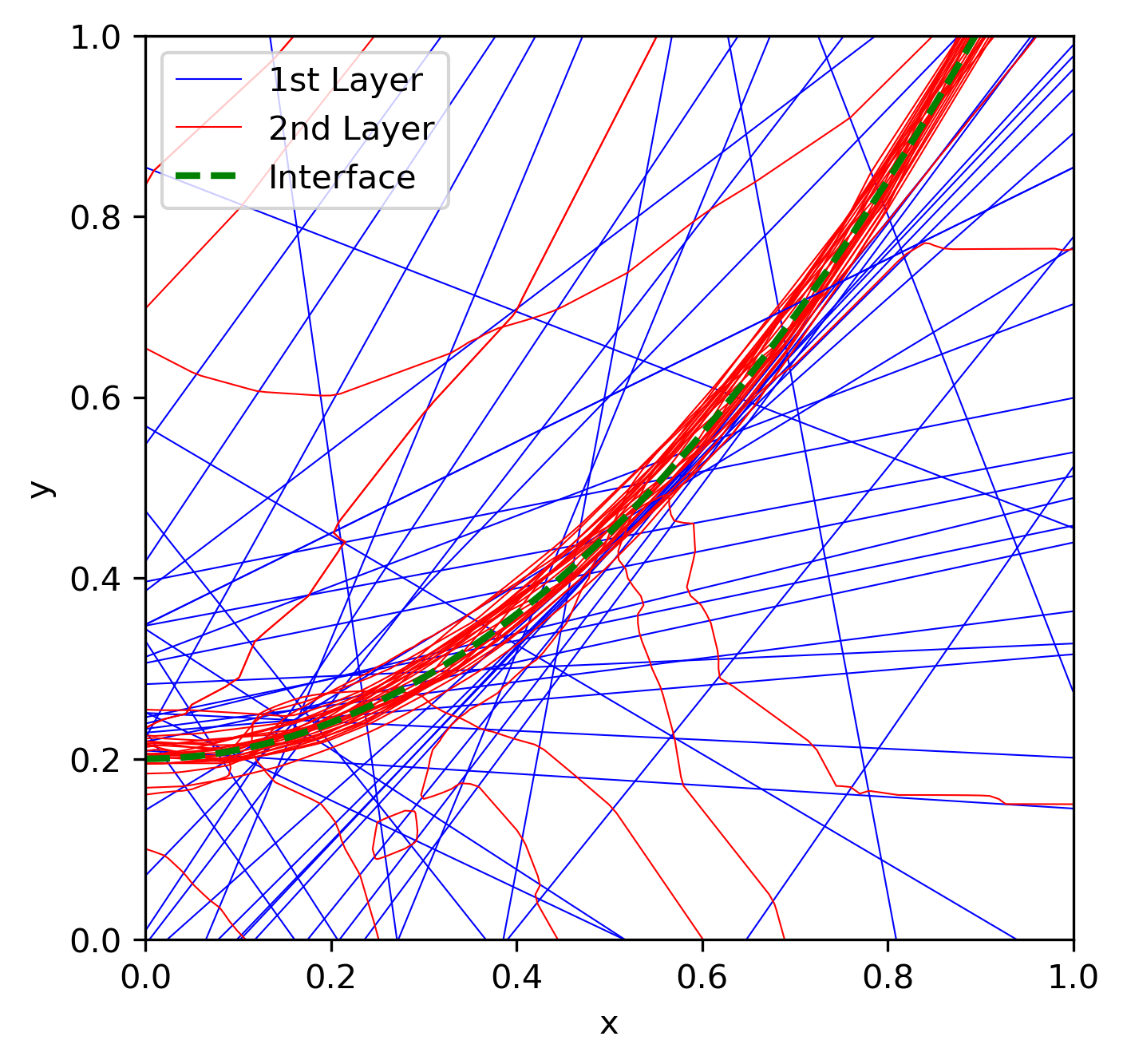}
\end{minipage}%
}%
\caption{Approximation results of the problem in \Cref{test5}}
\end{figure}

The LSNN method was implemented with 300000 iterations for 2--60--60--1 ReLU NN functions. We report the numerical results in \cref{test5 figure,test5 table}. Since the advective velocity field is a variable field, we increased the size of the neural network and $\rho=h/15$ was used to compute the finite difference quotient in \cref{finite_diff} to take values in one subdomain of the partition. Although theoretical analysis on the convergence of the method in the case of a smooth interface was not conducted, \cref{vertical5,comparison_exact5,comparison5,test5 table} show that the LSNN method is still capable of approximating the discontinuous solution with the curved interface accurately without oscillation. Finally, again, most of the second-layer breaking hyperplanes are along the interface (\cref{breaking5}) to approximate the discontinuous jump.

\begin{table}[htbp]\label{test5 table}
\caption{Relative errors of the problem in \Cref{test5}}
\centering
\begin{tabular}{|l|l|l|l|l|}
\hline
Network structure  &$\frac{\|u-{u}^{_N}_{_\cT}\|_0}{\|u\|_0}$ &$\frac{\vertiii{u-{u}^{_N}_{_\cT}}_{\bm\beta}}{\vertiii{u}_{\bm\beta}}$ & $\frac{\mathcal{L}^{1/2}({u}^{_N}_{_\cT},\bf f)}{\mathcal{L}^{1/2}({u}^{_N}_{_\cT},\bf 0)}$ & Parameters \\ \hline
2--60--60--1  & 0.046528 & 0.049423 & 0.019995   & 3901\\ \hline
\end{tabular}
\end{table}

\subsection{A problem with a constant advection velocity field ($d=3$)}\label{test6}
The last test problem is a three-dimensional problem defined on the domain $\Omega=(0,1)^3$, and approximation results are depicted on $z=0.505$. The advective velocity field is the constant field $\bm{\beta}(x,y,z) = (1,0,0),\,\, (x,y,z)\in\Omega$, and the inflow boundary of the problem is $\Gamma_{-}=\{(0,y,z):y,z\in(0,1)\}$.
The inflow boundary condition is given by
\begin{equation*}
g(x,y,z)=\left\{ \begin{array}{rl}
 1-4y,& (x,y,z)\in \Gamma^1_-\equiv\{(0,y,z)\in\Gamma_-: y<\frac{1}{2}\}, \\[2mm]
 5-4y, & (x,y,z)\in \Gamma^2_-=\Gamma_-\setminus \Gamma_-^1.
\end{array}\right.
\end{equation*} 
The exact solution of this test problem is 
\begin{equation}
u(x,y,z)=\left\{ \begin{array}{rl}
 1-4ye^{-x},& (x,y,z)\in \Omega_1\equiv\{(x,y,z)\in\Omega:y< \frac{1}{2}\}, \\[2mm]
 1+(4-4y)e^{-x}, & (x,y,z)\in \Omega_2=\Omega\setminus\Omega_1.
\end{array}\right.
\end{equation}

\begin{figure}[htbp]\label{test6 figure}
\centering
\subfigure[The interface\label{interface6}]{
\begin{minipage}[t]{0.4\linewidth}
\centering
\includegraphics[width=1.8in]{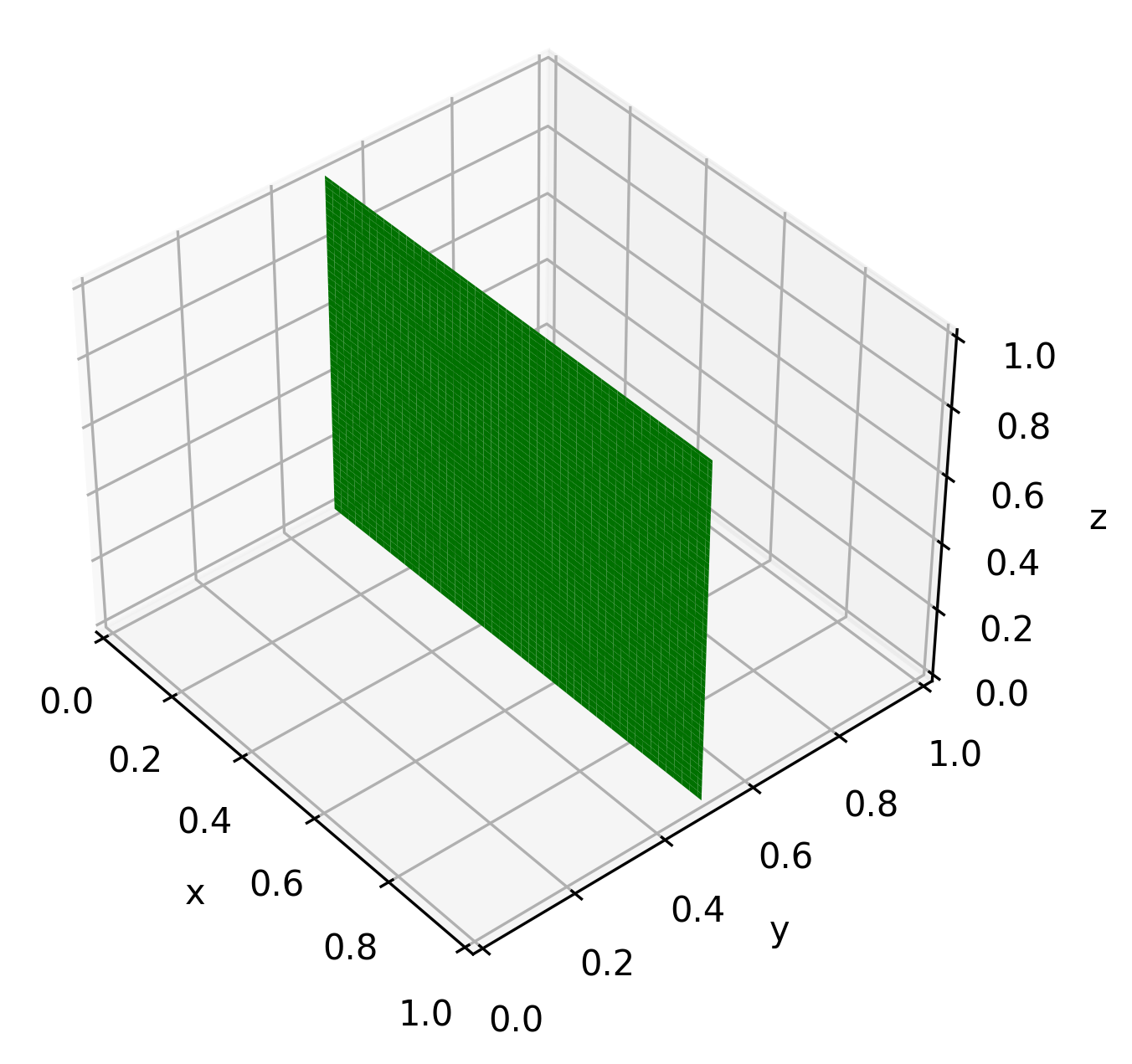}
\end{minipage}%
}%
\hspace{0.2in}
\subfigure[The trace of Figure \ref{comparison6} on $y=x$\label{vertical6}]{
\begin{minipage}[t]{0.4\linewidth}
\centering
\includegraphics[width=1.8in]{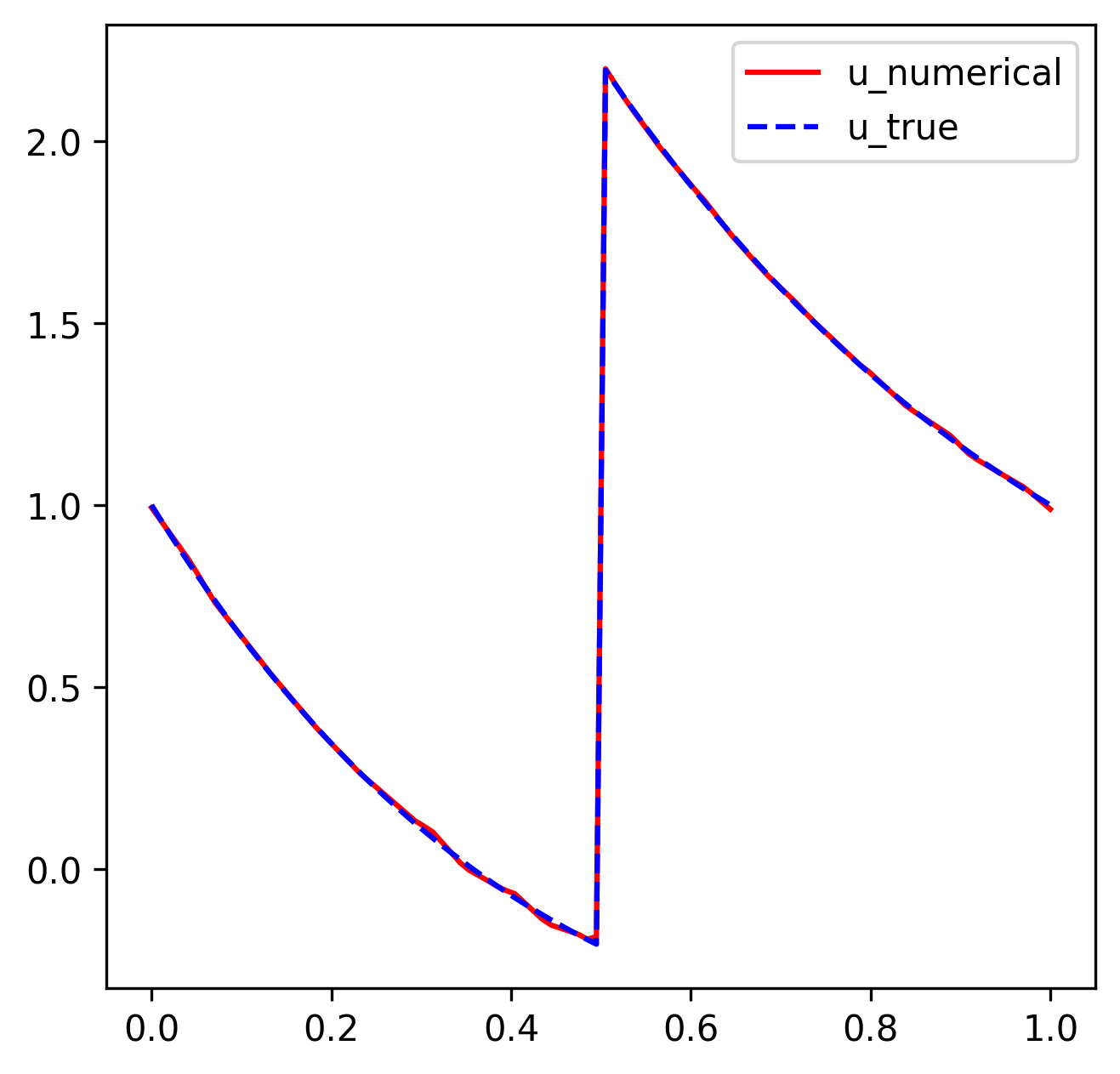}
\end{minipage}%
}%
\\
\subfigure[The exact solution\label{comparison_exact6}]{
\begin{minipage}[t]{0.4\linewidth}
\centering
\includegraphics[width=1.8in]{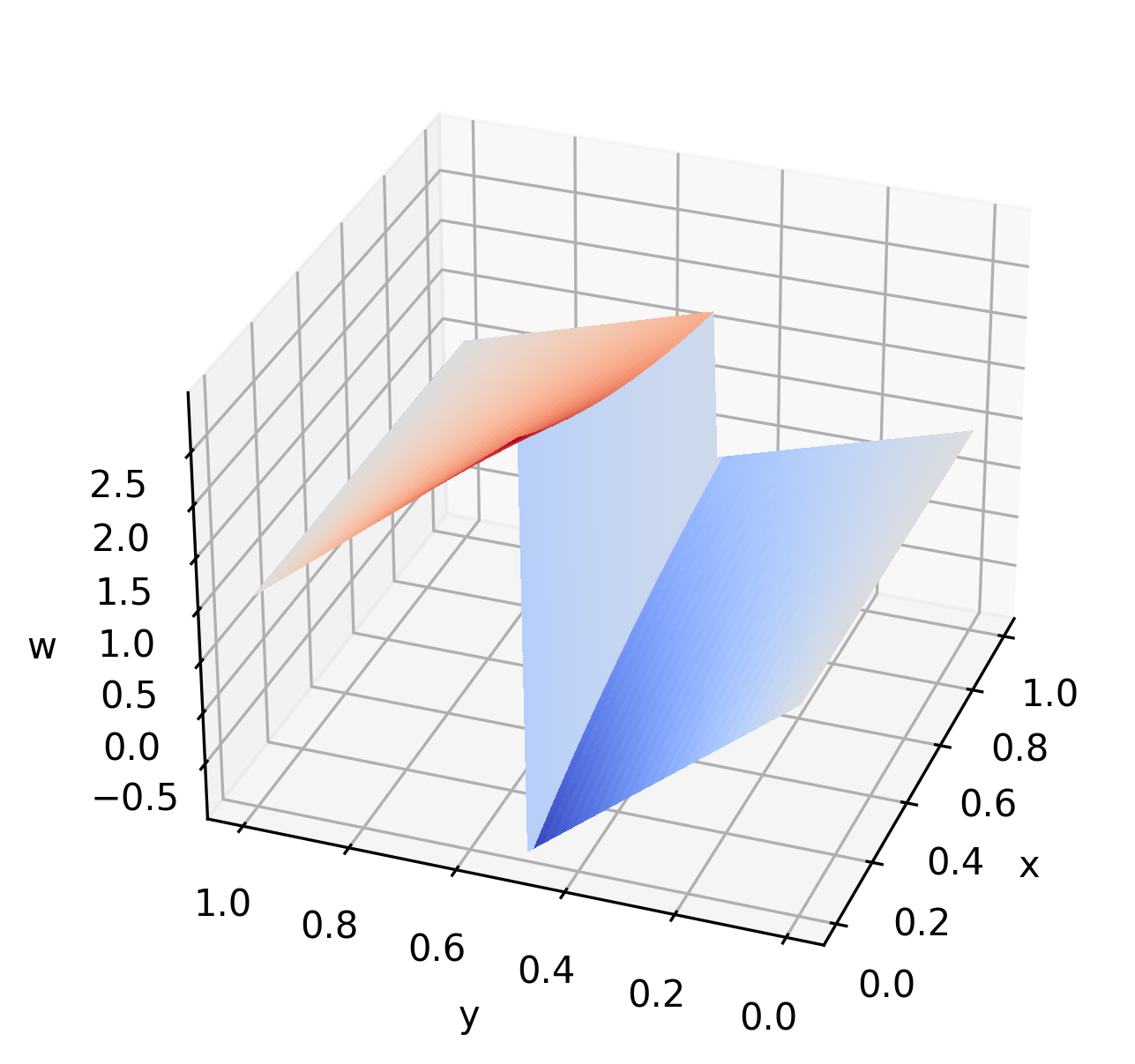}
\end{minipage}%
}%
\hspace{0.2in}
\subfigure[A 3--30--30--1 ReLU NN function approximation\label{comparison6}]{
\begin{minipage}[t]{0.4\linewidth}
\centering
\includegraphics[width=1.8in]{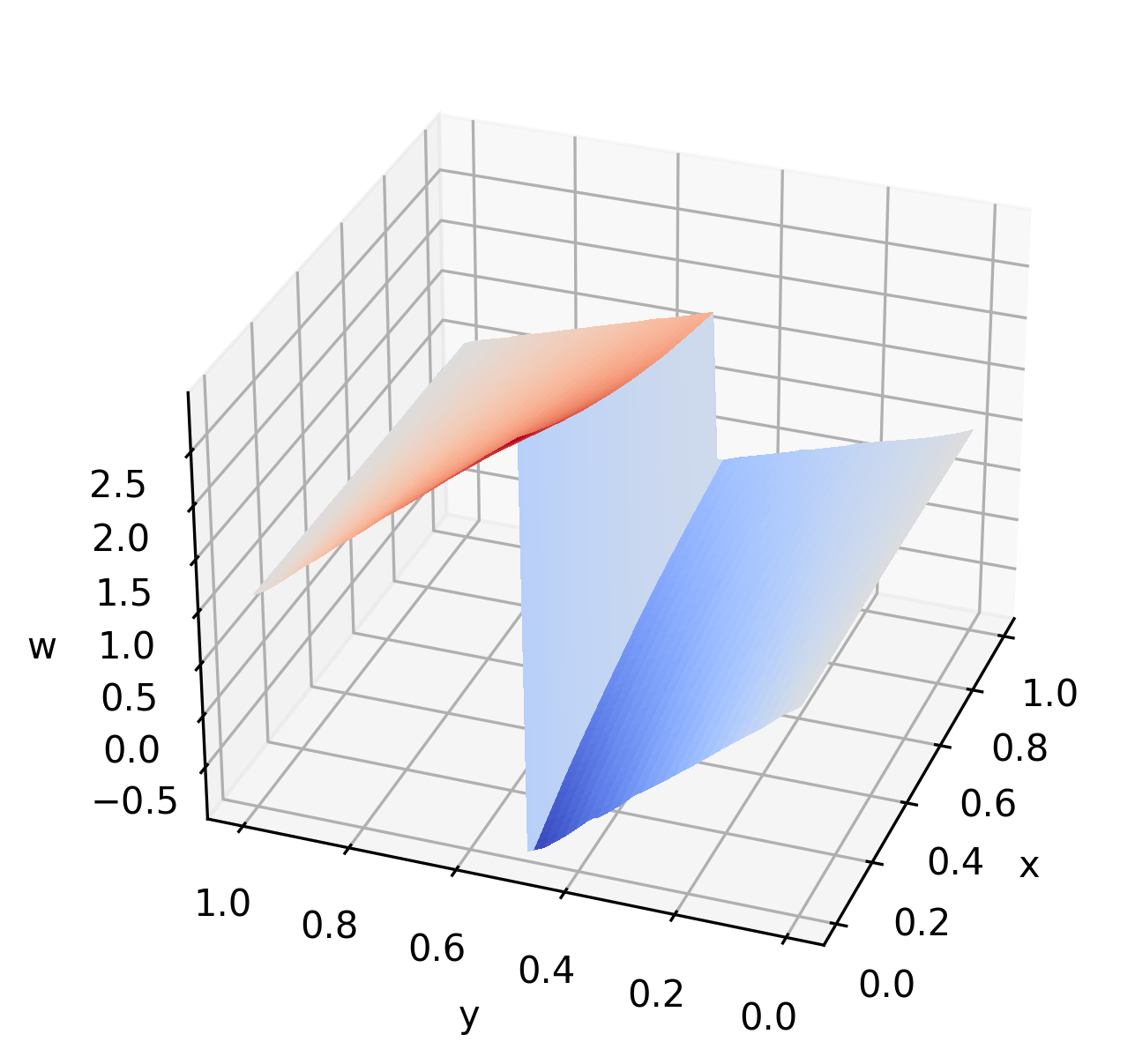}
\end{minipage}%
}%
\\
\subfigure[The breaking hyperplanes of the approximation in Figure \ref{comparison6}\label{breaking6}]{
\begin{minipage}[t]{0.4\linewidth}
\centering
\includegraphics[width=1.8in]{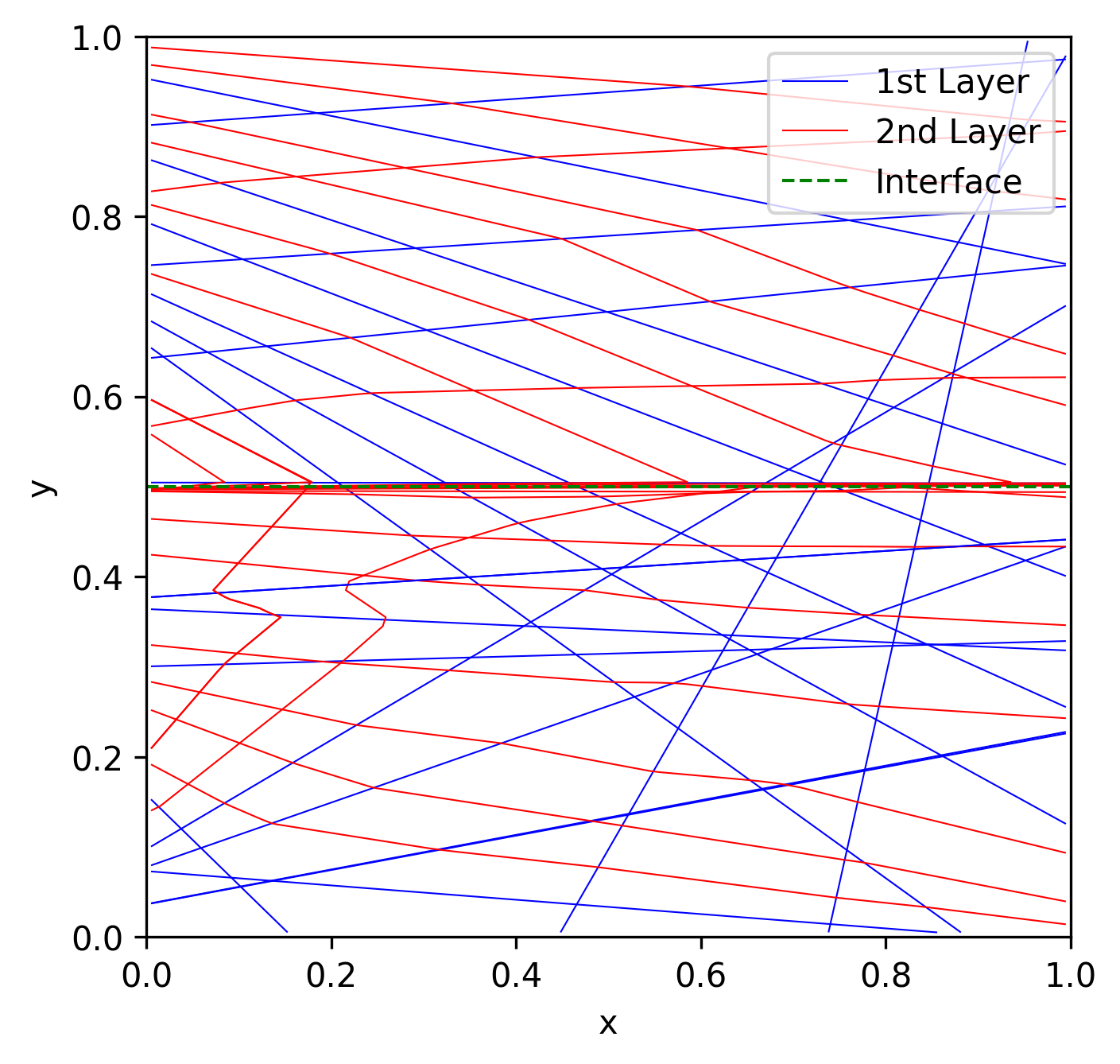}
\end{minipage}%
}%
\caption{Approximation results of the problem in \Cref{test6}}
\end{figure}

The LSNN method was implemented with 100000 iterations for 3--30--30--1 ReLU NN functions. Note that we still employ the two-hidden-layer NN since $\lceil \log_2(d+1)\rceil+1=3$ for $d=3$. We report the numerical results in \cref{test6 figure,test6 table}. The exact solution (\cref{comparison_exact6}) is piecewise smooth along the plane segment $y=0.5$ (\cref{interface6}) with a non-constant jump, and  was approximated accurately by the 3-layer NN (\cref{vertical6,comparison6,test6 table}). The behavior of the breaking hyperplanes (\cref{breaking6}) is similar to those of the previous examples around the discontinuous interface and on the subdomains.

\begin{table}[htbp]\label{test6 table}
\caption{Relative errors of the problem in \Cref{test6}}
\centering
\begin{tabular}{|l|l|l|l|l|}
\hline
Network structure  &$\frac{\|u-{u}^{_N}_{_\cT}\|_0}{\|u\|_0}$ &$\frac{\vertiii{u-{u}^{_N}_{_\cT}}_{\bm\beta}}{\vertiii{u}_{\bm\beta}}$ & $\frac{\mathcal{L}^{1/2}({u}^{_N}_{_\cT},\bf f)}{\mathcal{L}^{1/2}({u}^{_N}_{_\cT},\bf 0)}$ & Parameters \\ \hline
3--30--30--1  & 0.005082 & 0.027765 & 0.022519   & 1081\\ \hline
\end{tabular}
\end{table}

\section{Conclusion}\label{conclusion}
In this paper, we used the least-squares ReLU neural network (LSNN) method for solving linear advection-reaction
equation with discontinuous solution having non-constant jumps. The method, being mesh-free, requires no mesh to resolve the interfacial discontinuity. We proved theoretically that ReLU neural network (NN) functions with $\lceil \log_2(d+1)\rceil+1$-layer representations are capable of accurately approximating solutions with non-constant jumps along discontinuous interfaces that are not necessarily straight lines. Our theoretical findings were validated by multiple numerical examples with $d=2,3$ and various non-constant jumps and interface shapes, demonstrating accurate performance of the method.

The approximation of discontinuous functions by NNs is also encountered in classification tasks, and our results suggest that we can achieve accurate predictions with properly designed neural network architectures. However, as in \cite{cai2023least}, in this paper, we mainly focused on the depth of NNs. The approximation of discontinuous classification functions by NNs with fixed depth and width will be addressed in a forthcoming paper.

\bibliographystyle{siamplain}
\bibliography{references}

\section*{Appendix. The proof of Theorem \ref{thm:chi}}
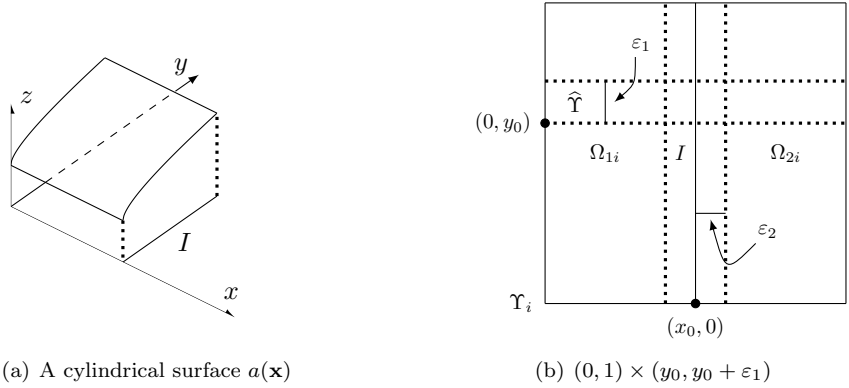
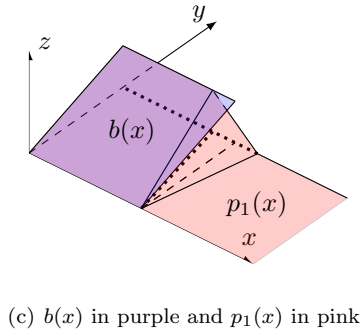
\begin{figure}[htbp]
\centering
\subfigure[A cylindrical surface $a(\mathbf{x})$\label{a(x)}]{
\begin{minipage}[t]{0.4\linewidth}
\centering
\begin{tikzpicture}
\begin{axis}[
axis lines=middle,
 view={40}{60},
 ticks=none,
  xmin=0,xmax=2.5,
  ymin=0,ymax=2.5,
  zmin=0,zmax=3,
  axis line style={draw=none}
]

\draw[-latex](0,0,0)--(2,0,0);
\node[] at (1.8,0.2,0) {$x$};

\draw[-latex](0,1.75,0)--(0,2,0);
\draw[dashed](0,0.37,0)--(0,1.75,0);
\draw[](0,0.37,0)--(0,0,0);
\node[] at (0.05,1.75,0.65) {$y$};

\draw[-latex](0,0,0)--(0,0,2.5);
\node[] at (0.1,0.05,2.7) {$z$};

\draw[samples=100,domain=0:1] plot(0,\x,{sqrt(\x)+1});
\draw[samples=100,domain=0:1] plot(1,\x,{sqrt(\x)+1});
\draw [] (0,0,1)-- (1,0,1);
\draw [] (0,1,2)-- (1,1,2);
\draw [] (1,0,0)-- (1,1,0);
\draw [dotted, very thick] (1,0,1)-- (1,0,0);
\draw [dotted, very thick] (1,1,0)-- (1,1,2);
\node[] at (1.3,0.3,0.4) {$I$};

\end{axis}
\end{tikzpicture}
\end{minipage}%
}%
\hspace{0.2in}
\subfigure[$(0,1)\times(y_0,y_0+\varepsilon_1)$\label{e1}]{
\begin{minipage}[t]{0.4\linewidth}
\centering
\begin{tikzpicture}[scale=0.8, transform shape]
    \draw [] (0,5)-- (5,5);
    \draw [] (0,5)-- (0,0);
    \draw [] (0,0)-- (5,0);
    \draw [] (5,0)-- (5,5);
    \draw [] (2.5,0)--(2.5,5);
    \draw [] (2.5,1.5)--(3,1.5);
    \draw [dotted, very thick] (2,0)--(2,5);
    \draw [dotted, very thick] (3,0)--(3,5);
    \draw [dotted, very thick] (0,3.7)--(5,3.7);
    \draw [dotted, very thick] (0,3)--(5,3);
    \draw [] (1,3)--(1,3.7);
    \draw[-latex] (3.5,1).. controls (3,0.5) ..(2.75,1.4);
    \draw[-latex] (1.5,4.1).. controls (1.5,3.65) ..(1.15,3.35);

    \node[] at (-0.4,0) {$\Upsilon_i$};
    \node[] at (.5,3.35) {$\widehat{\Upsilon}$};
    \node[] at (2.25,2.5) {$I$};
    \node[] at (1,2.5) {$\Omega_{1i}$};
    \node[] at (4,2.5) {$\Omega_{2i}$};
    \node[] at (2.5,-0.4) {$(x_0,0)$};
    \node[] at (-.7,3) {$(0,y_0)$};
    \filldraw[black] (2.5,0) circle (2pt);
    \filldraw[black] (0,3) circle (2pt);
    \node[] at (3.7,1.2) {$\varepsilon_2$};
    \node[] at (1.6,4.3) {$\varepsilon_1$};
    \end{tikzpicture}
\end{minipage}%
}%
\\
\subfigure[$b(x)$ in purple and $p_1(x)$ in pink\label{bp}]{
\begin{minipage}[t]{0.4\linewidth}
\centering
\begin{tikzpicture}
\begin{axis}[
axis lines=middle,
 view={40}{60},
 ticks=none,
  xmin=0,xmax=2.5,
  ymin=0,ymax=2.5,
  zmin=0,zmax=3,
  axis line style={draw=none}
]

\draw[-latex](0,0,0)--(2,0,0);
\node[] at (1.8,0.2,0) {$x$};

\draw[-latex](0,1.37,0)--(0,2,0);
\draw[dashed](0,0,0)--(0,1.37,0);
\node[] at (0.05,1.75,0.65) {$y$};

\draw[-latex](0,0,0)--(0,0,2.5);
\node[] at (0.1,0.05,2.7) {$z$};

\draw [] (0,0,0)-- (0,1,1);
\draw [] (0,1,1)-- (1,1,1);
\draw [] (1,0,0)-- (1,1,1);

\draw [] (0.935,1,0.663)-- (1.2,1,0);
\draw [dashed] (0.935,1,0.663)-- (0.8,1,1);
\draw [] (1,0,0)-- (1.2,1,0);
\draw [] (0.8,1,1)-- (1,0,0);
\draw [] (1.2,1,0)-- (2,1,0);

\draw [dashed] (1,0,0)-- (1,1,0);
\draw [dashed] (1,0,0)-- (1,1,0);
\draw [dotted, very thick] (1.2,1,0)-- (0,1,0);
\draw [dotted, very thick] (0.8,1,0)-- (1,0,0);

\draw[fill=red,opacity=0.2](1,0,0)--(1.2,1,0)--(2,1,0)--(2,0,0)--cycle;
\draw[fill=red,opacity=0.2](1,0,0)--(1.2,1,0)--(0.8,1,1)--cycle;
\draw[fill=red,opacity=0.2](1,0,0)--(0.8,1,1)--(0,1,1)--(0,0,0)--cycle;
\draw[fill=blue,opacity=0.2](1,0,0)--(1,1,1)--(0,1,1)--(0,0,0)--cycle;
\node[] at (0.5,0.5,0.4) {$b(x)$};
\node[] at (1.7,0.4,0.4) {$p_1(x)$};

\end{axis}
\end{tikzpicture}
\end{minipage}%
}%
\caption{An illustration of the convergence analysis on one subdomain $\Upsilon_i$}
\end{figure}

In this section, we provide a proof of \cref{thm:chi} by constructing a CPWL function to approximate $\chi(\bx)$ in \cref{chi}. We note that $a(\bx)$ is generally a cylindrical surface and that the jump of $\chi(\bx)$ is non-constant with \[
\chi(x,0)=a(x,0)=\alpha_1-\alpha_2 \text{ for } x\in [0,x_0]
\] (see \cref{a(x)}). 

For a given $\varepsilon_1>0$, we take $\widehat{\Upsilon}=(0,1)\times(y_0,y_0+\varepsilon_1)$ (see \cref{e1}). Without loss of generality, we let $\alpha_1=1$, $\alpha_2=0$, and $\widehat{\Upsilon}=(0,1)\times(0,\varepsilon_1)$.

Hence, 
\begin{equation}
    \chi(\bx)=\chi_0(\bx)+\chi_1(\bx) \text{ on }\widehat{\Upsilon},
\end{equation}
where $\chi_0(\bx)$ is a step function and $\chi_1(\bx)$ vanishes on the inflow boundary given by
\[
 \chi_0(\bx)=\left\{\begin{array}{rl}
 1, & \bx \in \Omega_{1i}\cap \widehat{\Upsilon},\\[2mm]
 0, & \bx \in \Omega_{2i}\cap \widehat{\Upsilon},
 \end{array}
 \right. \quad\mbox{and}\quad
 \chi_1(\bx)=\left\{\begin{array}{rl}
 a(\bx)-1, & \bx \in \Omega_{1i}\cap \widehat{\Upsilon},\\[2mm]
 0, & \bx \in \Omega_{2i}\cap \widehat{\Upsilon}.
 \end{array}
 \right.
\]

\begin{lemma}\label{lem:chi1}
Let
\[b(\bx)=\left\{\begin{array}{rl}
 \bm{b}\cdot (\bx-(x_0,0)), & \bx \in \Omega_{1i}\cap \widehat{\Upsilon},\\[2mm]
 0, & \bx \in \Omega_{2i}\cap \widehat{\Upsilon},
 \end{array}
 \right.\]
where $\bm{b}=(0,d)^T$ is a constant vector, and let $p_1(\bx)$ be the two-layer neural network function on $\widehat{\Upsilon}$ defined by
\[
p_1(\bx)=-c\,\sigma(\bw_1\cdot\bx+x_0)+c\,\sigma(\bw_2\cdot\bx+x_0)
\]
with the weights and coefficient
\begin{eqnarray*}
\bw_1=\begin{pmatrix}
    -1\\ -\varepsilon_2
\end{pmatrix}, \quad \bw_2=\begin{pmatrix}
    -1\\ \varepsilon_2
\end{pmatrix},\quad 
c=\dfrac{d}
{2\varepsilon_2}
\end{eqnarray*}
(see \cref{bp}). Then we have on $\widehat{\Upsilon}$,
 \begin{equation}\label{chi-est}
     \vertiii{b - p_1}_{\bm\beta}
     =\left(\|b - p_1\|^2_{0,\widehat{\Upsilon}} + \|b_{{\bm\beta}} - p_{1{\bm\beta}}\|^2_{0,\widehat{\Upsilon}}\right)^{1/2} 
     \le \sqrt{\frac{\varepsilon_1^3}{24}+\frac{B^2}{4}}\,|d|\, \sqrt{\varepsilon_1\varepsilon_2},
 \end{equation}
where we assume $|v_2(\mathbf{x})|\le B$ ($\bm\beta(\mathbf{x})=(0,v_2(\mathbf{x}))$).
\end{lemma}

\begin{proof}
Let us denote
\[
\widehat{\Upsilon}_{\varepsilon_2}\equiv\widehat{\Upsilon}_{1,\varepsilon_2}\cup \widehat{\Upsilon}_{2,\varepsilon_2}
\equiv \{\bx\in \Omega_{1i}\cap \widehat{\Upsilon} :\, \bw_1\cdot\bx+x_0<0\}\cup
\{\bx\in\Omega_{2i}\cap \widehat{\Upsilon} :\, \bw_2\cdot\bx+x_0>0\}.
\]
Then we have 
\[
b(\bx) -p_1(\bx) =\left\{\begin{array}{rl}
    -c(\bw_1\cdot\bx+x_0) & \bx\in\widehat{\Upsilon}_{1,\varepsilon_2},  \\[2mm]
    -c(\bw_2\cdot\bx+x_0) & \bx\in\widehat{\Upsilon}_{2,\varepsilon_2},  \\[2mm]
    0,  & \bx\in \widehat{\Upsilon}\setminus \widehat{\Upsilon}_{\varepsilon_2}.
\end{array}\right.
\]
Calculating the double integral,
\begin{equation}\label{chi-est1}
\|b - p_1\|^2_{0,\widehat{\Upsilon}}=\|b - p_1\|^2_{0,\widehat{\Upsilon}_{1,\varepsilon_2}}+\|b - p_1\|^2_{0,\widehat{\Upsilon}_{2,\varepsilon_2}}=\frac{d^2}{24}\varepsilon_1^4\varepsilon_2,
\end{equation}
and using the directional derivative,
\begin{equation}\label{chi-est2}
\|b_{{\bm\beta}} - p_{1{\bm\beta}}\|^2_{0,\widehat{\Upsilon}}=\int_{\widehat{\Upsilon}_{1,\varepsilon_2}}(c\bw_1\cdot\bm{\beta})^2\,d\bx+\int_{\widehat{\Upsilon}_{2,\varepsilon_2}}(c\bw_2\cdot\bm{\beta})^2\,d\bx\le\frac{(dB)^2}{4}\varepsilon_1\varepsilon_2.
\end{equation}
Now \cref{chi-est} follows from \cref{chi-est1,chi-est2}.
\end{proof}

\begin{lemma}\label{lem:chi2}
Let $\widehat{\Upsilon}$, $I$, $b(\mathbf{x})$, $p_1(\mathbf{x})$, and $\bm{\beta}(\mathbf{x})$ be as in {\em \cref{lem:chi1}} with $d=\chi_1(0,\varepsilon_1)/\varepsilon_1$, and let $p_0(\bx)$ be the two-layer neural network function on $\widehat{\Upsilon}$ defined by
\begin{equation}\label{p0}
p_0(\mathbf{x})=\frac{1}{2\varepsilon_2}\left(\sigma(x-x_0+\varepsilon_2)-\sigma(x-x_0-\varepsilon_2)\right).  
\end{equation}
Then we have on $\widehat{\Upsilon}$,
\begin{equation}\label{lem:chi2 ineq}
     \vertiii{\chi - (p_0+p_1)}_{\bm\beta}\le C_1\sqrt{\varepsilon_1\varepsilon_2}+C_2\sqrt{\varepsilon_1},
\end{equation}
where $C_2$ is given by the square root of
\begin{multline}
    \sup\{(u_1(\mathbf{x})-u_2(\mathbf{x})-1-b(\mathbf{x}))^2:\mathbf{x}\in\Omega_{1i}\cap \widehat{\Upsilon}\}x_0\\+\sup\{2\gamma(\mathbf{x})^2(u_2(\mathbf{x})-u_1(\mathbf{x}))^2+2(dB)^2:\mathbf{x}\in \Omega_{1i}\cap \widehat{\Upsilon}\}x_0.
\end{multline}
\end{lemma}

\begin{proof}
From the triangle inequality,
\begin{equation}
\vertiii{\chi - (p_0+p_1)}_{\bm\beta}=\vertiii{\chi_0+\chi_1 - (p_0+p_1)}_{\bm\beta}\le \vertiii{\chi_0 - p_0}_{\bm\beta}+\vertiii{\chi_1 - p_1}_{\bm\beta}.
\end{equation}
Since $\|\chi_{0\bm\beta} - p_{0\bm\beta}\|=0$, calculating the double integral,
\begin{equation}
\begin{split}
\vertiii{\chi_0 - p_0}_{\bm\beta}
     &=\left(\|\chi_0 - p_0\|^2_{0,\widehat{\Upsilon}} + \|\chi_{0\bm\beta} - p_{0\bm\beta}\|^2_{0,\widehat{\Upsilon}}\right)^{1/2}\\
     &=\|\chi_0 - p_0\|_{0,\widehat{\Upsilon}}\\
     &=\frac{1}{\sqrt{6}}\, \sqrt{\varepsilon_1\varepsilon_2}.
\end{split}
\end{equation}
Next, again, by the triangle inequality,
\[\vertiii{\chi_1 - p_1}_{\bm\beta}\le \vertiii{\chi_1 - b}_{\bm\beta}+\vertiii{b - p_1}_{\bm\beta}.\]
By \cref{lem:chi1}
\[\vertiii{b - p_1}_{\bm\beta}\le \sqrt{\frac{\varepsilon_1^3}{24}+\frac{B^2}{4}}\,|d|\, \sqrt{\varepsilon_1\varepsilon_2}.\]
To bound $\vertiii{\chi_1 - b}_{\bm\beta}$, we recall the definition of the graph norm,
\[\vertiii{\chi_1 - b}_{\bm\beta}=\left(\|\chi_1 - b\|^2_{0,\widehat{\Upsilon}} + \|\chi_{1\bm\beta} - b_{\bm\beta}\|^2_{0,\widehat{\Upsilon}}\right)^{1/2}.\]
First we have
\begin{equation*}
\begin{split}
 \|\chi_1 - b\|^2_{0,\widehat{\Upsilon}}&=\|\chi_1 - b\|^2_{0,\Omega_{1i}\cap \widehat{\Upsilon}}\\
 &\le \sup\{(\chi_1(\mathbf{x})-b(\mathbf{x}))^2:\mathbf{x}\in\Omega_{1i}\cap \widehat{\Upsilon}\} \varepsilon_1x_0\\
 &= \sup\{(u_1(\mathbf{x})-u_2(\mathbf{x})-1-b(\mathbf{x}))^2:\mathbf{x}\in\Omega_{1i}\cap \widehat{\Upsilon}\} \varepsilon_1x_0.   
\end{split}
\end{equation*}
Next, observing $b_{\bm\beta}=dv_2$ and from \cref{pde},
\[\chi_{1\bm\beta}=(u_1-u_2-1)_{\bm\beta}=(u_1-u_2)_{\bm\beta}=\gamma(u_2-u_1),\]
we have
\begin{equation*}
\begin{split}
\|\chi_{1\bm\beta} - b_{\bm\beta}\|^2_{0,\widehat{\Upsilon}}&=\|\chi_{1\bm\beta} - b_{\bm\beta}\|^2_{0,\Omega_{1i}\cap \widehat{\Upsilon}}\\ 
&=\|\gamma(u_2-u_1)-dv_2\|^2_{0,\Omega_{1i}\cap \widehat{\Upsilon}}\\
&\le\left\lVert\sqrt{2\gamma^2(u_2-u_1)^2+2(dv_2)^2}\right\rVert^2_{0,\Omega_{1i}\cap \widehat{\Upsilon}}\\
&\le \sup\{2\gamma(\mathbf{x})^2(u_2(\mathbf{x})-u_1(\mathbf{x}))^2+2(dB)^2:\mathbf{x}\in \Omega_{1i}\cap \widehat{\Upsilon}\}\varepsilon_1x_0.
\end{split}
\end{equation*}
Now \cref{lem:chi2 ineq} follows from combining the above inequalities.
\end{proof}

Given $\varepsilon_3>0$, let us choose $\varepsilon_1=1/m$ such that
\begin{multline}
\sup\{(u_1(\bx)-u_2(\bx)-\chi(0,j/m)-b_j(\mathbf{x}))^2:\mathbf{x}\in\Omega_{1i}\cap \widehat{\Upsilon}_j\},\\
\sup\{2\gamma(\mathbf{x})^2(u_2(\mathbf{x})-u_1(\mathbf{x}))^2+2(d_jB)^2:\mathbf{x}\in\Omega_{1i}\cap \widehat{\Upsilon}_j\}<\varepsilon_3,
\end{multline}
where $\widehat{\Upsilon}_j=(0,1)\times(j/m,(j+1)/m)$ for $j=0,\ldots,m-1$, and $b_j(\mathbf{x})$ and $d_j$ are as in \cref{lem:chi2}. Then we define on each $\widehat{\Upsilon}_j$, $p_{0j}(\mathbf{x}),p_{1j}(\mathbf{x})$ as in \cref{lem:chi2}, and construct the CPWL function $p_i(\mathbf{x})$ on $\Upsilon_i$ defined by
\[p_i(\mathbf{x})=p_{0j}(\mathbf{x})+p_{1j}(\mathbf{x}),\ \mathbf{x}\in\widehat{\Upsilon}_j.\]

\begin{proof}[Proof of \cref{thm:chi}]
By \cref{lem:chi2} and the given condition,
\begin{equation}
\begin{split}
\vertiii{\chi-p_i}_{\bm\beta}&=\left(\sum_{j=0}^{m-1}\vertiii{\chi-(p_{0j}+p_{1j})}^2_{\bm\beta}\right)^{1/2}\\  &\le \left(\sum_{j=0}^{m-1} \left(D_{1j}\sqrt{\varepsilon_1\varepsilon_2}+D_{2j}\sqrt{\varepsilon_1\varepsilon_3}\right)^2\right)^{1/2}\\
&\le\left(m\max_{0\le j\le m-1}\left(D_{1j}\sqrt{\varepsilon_1\varepsilon_2}+D_{2j}\sqrt{\varepsilon_1\varepsilon_3}\right)^2\right)^{1/2}\\
&= \sqrt{m}\max_{0\le j\le m-1}\left(D_{1j}\sqrt{\varepsilon_1\varepsilon_2}+D_{2j}\sqrt{\varepsilon_1\varepsilon_3}\right)\\
&=\max_{0\le j\le m-1}\left(D_{1j}\sqrt{\varepsilon_2}+D_{2j}\sqrt{\varepsilon_3}\right),
\end{split}
\end{equation}
where for the first identity, each norm on the right-hand side is taken over $\widehat{\Upsilon}_j$. Now $D_1=D_{1k}$ and $D_2=D_{2k}$ for some $0\le k\le m-1$.
\end{proof}
\end{document}